\documentclass[a4paper]{amsart}

\usepackage[english]{babel}

\usepackage{amssymb}
\usepackage{amsthm}
\usepackage{mathtools}
\usepackage{dsfont}
\usepackage{enumitem}
\usepackage{aliascnt}
\usepackage{enumitem}
\usepackage{amsmath}
\usepackage{longtable}
\usepackage{booktabs}
\usepackage{tabularx}
\usepackage{changepage}

\usepackage[T1]{fontenc}
\usepackage{lmodern}
\usepackage{microtype}

\usepackage[a4paper, inner=3cm,outer=3cm,top=3.3cm,bottom=3cm]{geometry}

\usepackage{xspace}

\usepackage[autostyle]{csquotes}

\usepackage{tikz}
\usetikzlibrary{matrix,arrows,patterns,decorations.pathreplacing}
\usepackage{pgfplots}

\usepackage[colorlinks, bookmarksdepth=3,citecolor=blue]{hyperref}

\numberwithin{equation}{section}

\makeatletter
\@ifpackageloaded{hyperref}{%
\DeclareRobustCommand{\SkipTocEntry}[5]{}}{%
\DeclareRobustCommand{\SkipTocEntry}[4]{}}
\makeatother

\newtheoremstyle{mythmstyle}
{6pt} 
{6pt} 
{\itshape}
{}
{\scshape}
{.}
{3pt}
{}
\theoremstyle{mythmstyle}

\newtheorem{thm}{Theorem}[section]
\newaliascnt{prop}{thm}
\newtheorem{prop}[prop]{Proposition}
\aliascntresetthe{prop}
\newaliascnt{cor}{thm}
\newtheorem{cor}[cor]{Corollary}
\aliascntresetthe{cor}
\newaliascnt{lem}{thm}
\newtheorem{lem}[lem]{Lemma}
\aliascntresetthe{lem}
\newaliascnt{conj}{thm}

\aliascntresetthe{conj}

\newaliascnt{enumprob}{thm}
\newtheorem{enumprob}[enumprob]{Enumeration Problem}
\aliascntresetthe{enumprob}

\newaliascnt{classprob}{thm}
\newtheorem{classprob}[classprob]{Classification Problem}
\aliascntresetthe{classprob}

\newtheoremstyle{mydefnstyle}
{6pt} 
{6pt} 
{\upshape}
{}
{\scshape}
{.}
{3pt}
{}
\theoremstyle{mydefnstyle}

\newtheorem{defn}[thm]{Definition}
\newtheorem{ex}[thm]{Example}
\newaliascnt{rem}{thm}
\newtheorem{rem}[rem]{Remark}
\aliascntresetthe{rem}
\newaliascnt{algo}{thm}
\newtheorem{algo}[rem]{Algorithm}
\aliascntresetthe{algo}

\newcommand{\cP}{\mathcal{P}}
\newcommand{\vertset}{\operatorname{vert}}
\newcommand{\lexmin}{\operatorname{lexmin}}

\DeclareMathOperator{\conv}{conv}

\DeclareMathOperator{\Id}{Id}

\newcommand{\N}{{\mathds{N}}}

\newcommand{\R}{{\mathds{R}}}
\newcommand{\Z}{{\mathds{Z}}}

\newcommand{\C}{\mathds{C}}

\newcommand{\rleft}{\mathopen{}\mathclose\bgroup\left}
\newcommand{\rright}{\aftergroup\egroup\right}

\newcommand{\setcond}[2]{\left\{ #1 \,:\, #2 \right\}}

\newcommand{\GL}{\operatorname{GL}}
\newcommand{\Aff}{\operatorname{Aff}}
\newcommand{\GG}{\mathbb{G}}
\newcommand{\vol}{\operatorname{vol}}
\newcommand{\Cay}{\operatorname{Cay}}

\DeclareMathOperator{\aff}{aff}

\newcommand{\set}[1]{\rleft\{ {#1} \rright\}}

\DeclareMathOperator{\Vol}{Vol}

\newcommand{\cal}[1]{\mathcal{#1}}

\newcommand{\mS}{\mathbb S}

\newcommand{\D}{\Delta}

\newcommand{\cF}{\cal F}

\newcommand{\cS}{\cal S}

\newcommand{\supp}{\operatorname{supp}}
\newcommand{\la}{\langle}
\newcommand{\ra}{\rangle}

\newcommand{\rs}[1]{Section~\ref{S:#1}}

\newcommand{\rp}[1]{Proposition~\ref{P:#1}}

\newcommand{\re}[1]{(\ref{eq:#1})}

\newcommand{\rt}[1] {Theorem~\ref{T:#1}}
\newcommand{\rd}[1]{Definition~\ref{D:#1}}

\newcommand{\ralgo}[1]{Algorithm~\ref{algo:#1}}

\newcommand\blfootnote[1]{%
  \begingroup
  \renewcommand\thefootnote{}\footnote{#1}%
  \addtocounter{footnote}{-1}%
  \endgroup
}

\title[Classification of triples of lattice polytopes with a given mixed volume]{Classification of triples of lattice polytopes\\ with a given mixed volume}
\author{Gennadiy~Averkov}
\address[Gennadiy~Averkov]{Fakult\"at 1, BTU Cottbus-Senftenberg, Platz der Deutschen Einheit 1, 03046 Cottbus, Germany}
\email{averkov@b-tu.de}
\author{Christopher~Borger}
\address[Christopher~Borger]{Faculty of Mathematics, Otto-von-Guericke-Universit\"at Magdeburg, Universit\"atsplatz 2, 39106 Magdeburg, Germany}
\email{christopher.borger@ovgu.de}
\author{Ivan~Soprunov}
\address[Ivan~Soprunov]{Department of Mathematics and Statistics, Cleveland State University,  2121 Euclid Ave, Cleveland, Ohio, 44115 USA}
\email{i.soprunov@csuohio.edu}

\usepackage{lipsum}

\begin{document}
\selectlanguage{english}

\begin{abstract}
We present an algorithm for the classification of
triples of lattice polytopes with a given mixed
volume $m$ in dimension $3$. 
It is known that the classification can be reduced 
to the enumeration of so-called irreducible triples, the number of 
which is finite for fixed $m$. Following this algorithm, we 
enumerate all irreducible triples of normalized mixed volume up to $4$
that are inclusion-maximal. This produces a classification of generic
trivariate sparse polynomial systems with up to 4 solutions in the complex torus, up to monomial changes of variables.
By a recent result of Esterov, this leads to a description of all generic trivariate sparse polynomial systems that are solvable by radicals. 
\end{abstract}



\maketitle


\section{Introduction} 

\subsection{Formulation of the problem and previous results}

In the last decade, there has been an increased interest in the algorithmic theory of lattice polytopes, which is motivated by applications in algebra, algebraic geometry, combinatorics, and optimization (see, for example, \cite{AKW17,averkov2011maximal,balletti2019families,blanco2016finiteness,blanco2016lattice,blanco2018enumeration,LorenzPaffenholz15,castryck2012moving,iglesias2019complete,Kasprzyk10, KreuzerSkarke1998,NillObro10}). So far, a special emphasis has been put on computer-assisted enumeration results, which are important from different perspectives. On the one hand, one can carry out enumeration to gather concrete data and then make and verify hypotheses based on these data. On the other hand, proving results on lattice polytopes frequently requires to handle special cases, which are only finitely many but are hard to determine without algorithmic assistance. Most notably, some
structural classification results for lattice polytopes hold up to finitely many exceptional situations and, in this case, one is interested in describing the exceptional situations by means of enumeration in order to accomplish a classification. Thus, the structural and algorithmic theory are intertwined and constantly influence each other. Evaluation of new data established using a computer-assisted search leads to new theoretical questions, while new theoretical results (in particular, finiteness results) suggest new enumeration tasks. 

Our point of departure is the classical Bernstein--Khovanskii--Kouchnirenko theorem that determines the number of solutions of a generic system of Laurent polynomial equations:

\begin{thm}[\cite{Bernstein75}]
\label{thm:BKK} Let $A_1,\dots, A_d$ be non-empty finite subsets of $\Z^d$ and 
consider systems of $d$ equations $f_1=\cdots = f_d =0$ where each $f_i$ is a Laurent polynomial of the form 
\[
	f_i(x_1,\ldots,x_d) = \sum_{(a_1,\ldots,a_d) \in A_i} c_{i,(a_1,\ldots,a_d)} x_1^{a_1} \cdots x_d^{a_d}.
\] 
Then, for a generic choice of the coefficients $c_{i,a} \in \C$ $(1 \le i \le d$, $a \in A_i)$, the number of solutions of the system $f_1= \cdots = f_d = 0$ in the complex torus $(\C \setminus \{0\})^d$ is equal to the normalized mixed volume $V(P_1,\ldots,P_d)$ of the lattice polytopes $P_1 =\conv(A_1),\ldots, P_d = \conv(A_d)$.
\end{thm}

The subset $A_i$ is called the \emph{support} of $f_i$ and the polytope $P_i=\conv(A_i)$ is the \emph{Newton polytope} of $f_i$. Systems with 
generic coefficients and prescribed supports are often called {\it sparse} as they may have much fewer monomials than generic systems with prescribed degrees of the polynomials. By  \emph{generic choice}  of the coefficients  we mean that the vector of all coefficients of the system is chosen outside an appropriate proper algebraic subset of $\C^{|A_1| + \cdots + |A_d|}$. There is extensive literature on sparse systems covering both computational and theoretical points of view. For example, the reader may consult \cite[Chapter 7]{CLO04} for a list of references. 

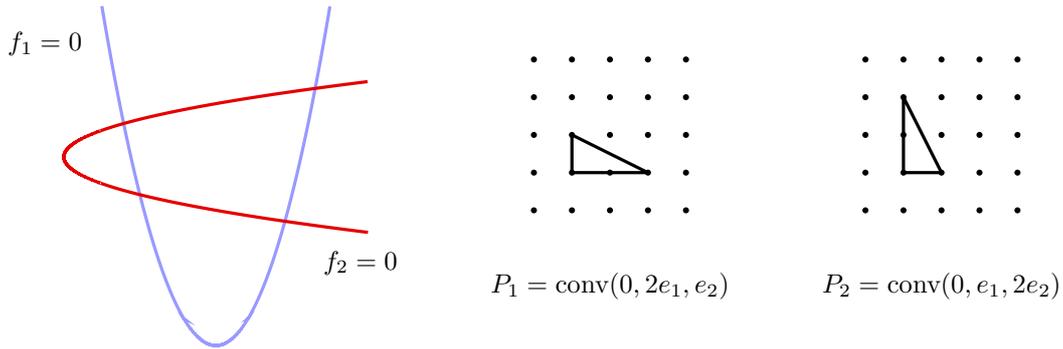
\begin{figure}[h]
\begin{tikzpicture}[scale=0.5,baseline=-20mm]
\draw [rounded corners=.3mm, very thick,blue!40!white,smooth,samples =100,domain=-3:3] plot(\x,{\x*\x-5});
\node at (-4.5,3){$f_1=0$};
\node at (3.8,-2.8){$f_2=0$};
\draw [rounded corners=.3mm, very thick,red!90!black,smooth,samples =100,domain=-2:2] plot({2*\x*\x-4},\x);
\end{tikzpicture}
\hspace{5ex}
\begin{tikzpicture}[scale=0.5]
	\draw[rounded corners=.3mm, very thick] (0,0) -- (2,0) -- (0,1) -- cycle;
	\node at (1, -3) {$P_1=\conv(0,2 e_1, e_2)$};
	\foreach \x in {-1,...,3}
	\foreach \y in {-1,...,3}
	{
		\fill (\x,\y) circle (0.08);		
	}
\end{tikzpicture}
\hspace{5ex}
\begin{tikzpicture}[scale=0.5]
	\draw[rounded corners=.3mm, very thick] (0,0) -- (1,0) -- (0,2) -- cycle;
	\node at (1, -3) {$P_2=\conv(0, e_1, 2 e_2)$};
	\foreach \x in {-1,...,3}
\foreach \y in {-1,...,3}
{
	\fill (\x,\y) circle (0.08);		
}
\end{tikzpicture}
\caption{The system $f_1=f_2=0$, where $f_1=0$ is a generic vertical parabola and $f_2=0$ is a generic horizontal parabola has $4$ solutions in $(\C \setminus \{0\})^2$, because the normalized mixed volume of $P_1$ and $P_2$ is $4$. \label{fig:ex:BKK}}
\end{figure}

\begin{ex}
	\label{ex:BKK}
	Consider general equations $f_1=0$ and $f_2=0$ of vertically and horizontally aligned parabolas given  by polynomials 
	\begin{align*}
		f_1(x,y) & = c_{1,(0,0)}+ c_{1,(1,0)} x +c_{1,(2,0)} x^2 + c_{1,(0,1)} y,
		\\ f_2(x,y) & = c_{2,(0,0)}+c_{2,(1,0)} x +  c_{2,(0,1)} y +c_{2,(0,2)} y^2. 
	\end{align*}
	In this case the supports $A_1$ and $A_2$ are 
	\begin{align*}
		A_1 & = \{(0,0), (1,0), (2,0), (0,1)\},
		&  A_2 & = \{(0,0), (1,0), (0,1),(0,2)\}
	\end{align*} 
	and the Newton polytopes  $P_1$ and $P_2$ are triangles 
	\begin{align*}
		P_1 & =\conv(A_1) = \conv((0,0),(2,0),(1,0)),
		& P_2 & = \conv(A_2) = \conv((0,0),(1,0),(0,2)).
	\end{align*} By Theorem~\ref{thm:BKK}, if the vector 
	\[
		(c_{1,(0,0)},c_{1,(1,0)}, c_{1,(2,0)},c_{1,(0,1)}, c_{2,(0,0)}, c_{2,(1,0)}, c_{2,(0,1)},c_{2,(0,2)}) \in \C^8
	\] 
	of all coefficients of the polynomials $f_1$ and $f_2$
	is generic, then the system $f_1=f_2=0$ has exactly $4$ solutions in $(\C \setminus \{0\})^2$, because the normalized mixed volume $V(P_1,P_2)$ of $P_1$ and $P_2$ equals $4$. The value $V(P_1,P_2)$ can be determined from the inclusion-exclusion type formula 
	\[
		V(P_1,P_2) = \frac{1}{2} ( \Vol(P_1 + P_2) - \Vol(P_1) - \Vol(P_2)).
	\]
	Here $\Vol(P)$ denotes the normalized 2-dimensional volume of a polytope $P$, which is twice the Euclidean area of $P$, 
	 and $P_1+ P_2$ is the Minkowski sum of the two triangles.  See also Fig.~\ref{fig:ex:BKK}.
\end{ex} 

In the notation of Theorem~\ref{thm:BKK}, we call $(A_1,\ldots,A_d)$ the (total) \emph{support} of the system $f_1=\cdots = f_d=0$. Theorem~\ref{thm:BKK} determines the number of solutions of a generic sparse system with a given support.  Recently, a reverse research direction was suggested by Esterov and Gusev \cite{EsterovGusev2015,EsterovGusev2016}: One can fix the number of solutions $m$ and try to classify all supports $(A_1,\ldots,A_d)$ of generic systems with exactly $m$ solutions. By Theorem~\ref{thm:BKK}, it suffices to use the convex hulls of the sets $A_1,\ldots,A_d$ in such a classification. This leads to the following problem.
\begin{classprob}
	\label{class:prob}
	Given $d, m \in \N$, describe all $d$-tuples $(P_1,\ldots,P_d)$ of lattice polytopes whose normalized mixed volume equals $m$.
\end{classprob}

The solution of this problem has been known only in the following cases: 
\begin{itemize}
	\item The trivial case $d=1$.
	\item The case $m=1$, for which a solution is provided by a result of Esterov and Gusev  \cite{EsterovGusev2015}.
	\item The case $d=2, m \le 4$, for which a solution is described by Esterov and Gusev in \cite{EsterovGusev2016}. 
\end{itemize}

Classification~Problem~\ref{class:prob} is of particular interest for $m \le 4$, since 
by a result of Esterov \cite{Esterov2018}, it allows to describe all generic system $f_1=\cdots = f_d=0$ as in Theorem~\ref{thm:BKK} which are solvable by radicals. 

The family of $d$-tuples of lattice polytopes with mixed volume $m$ is invariant under application of a common unimodular transformation to all polytopes of the tuple, independent translations of the polytopes by lattice vectors, and permutations of the polytopes of the $d$-tuple. We call $d$-tuples \emph{equivalent} if they coincide up to the above transformations.

Note that, even modulo the above equivalence relation, the set of such $d$-tuples is not finite. Nevertheless, the classification can be reduced to the enumeration of so-called irreducible tuples.   Recall that
a tuple $(P_1,\dots,P_d)$ of polytopes in $\R^d$ is called {\it irreducible} if the Minkowski sum
of any $k$ polytopes in the tuple is at least $(k+1)$-dimensional for every $1\leq k\leq d-1$. In  \cite{Esterov2018}
Esterov showed that the number of irreducible $d$-tuples of lattice polytopes is
finite up to equivalence for a fixed mixed volume and dimension
(see Theorem~\ref{T:finiteness} below). 

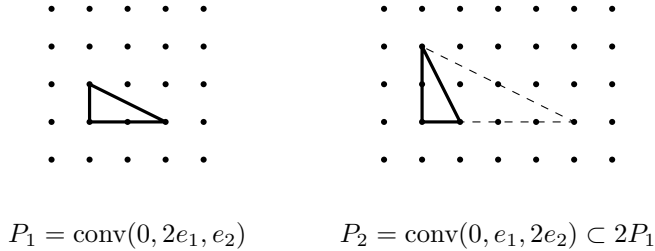
\begin{figure}[h]
\begin{tikzpicture}[scale=0.5]
\draw[rounded corners=.3mm, very thick] (0,0) -- (2,0) -- (0,1) -- cycle;
\node at (1, -3) {$P_1=\conv(0,2 e_1, e_2)$};
\foreach \x in {-1,...,3}
\foreach \y in {-1,...,3}
{
	\fill (\x,\y) circle (0.08);		
}
\end{tikzpicture}
\hspace{5ex}
\begin{tikzpicture}[scale=0.5]
\draw[rounded corners=.3mm, very thick] (0,0) -- (1,0) -- (0,2) -- cycle;
\draw[rounded corners=.3mm, dashed] (0,0) -- (4,0) -- (0,2) -- cycle;
\node at (2, -3) {$P_2=\conv(0, e_1, 2 e_2) \subset 2 P_1$};
\foreach \x in {-1,...,5}
\foreach \y in {-1,...,3}
{
	\fill (\x,\y) circle (0.08);		
}
\end{tikzpicture}
\caption{Embedding the pair $(P_1,P_2)$ from Example~\ref{ex:BKK} into a maximal pair by enlarging $P_2$\label{ex:BKK:enlarging} to $2 P_1$. Since $V(P_1, 2 P_1) = 2 V(P_1,P_1) = 2 \Vol(P_1) = 4$, the normalized mixed volume remains unchanged.}
\end{figure}

There is a natural partial ordering on the set of equivalence classes of $d$-tuples of lattice polytopes with a given mixed volume $m$, 
induced by inclusion. Maximal elements with respect to this order are called {\it maximal} $d$-tuples with a given mixed volume. Theorem~\ref{T:finiteness} implies that every irreducible $d$-tuple of lattice polytopes is contained in a maximal irreducible $d$-tuple of lattice polytopes. As a consequence, enumeration of irreducible $d$-tuples amounts to enumeration of maximal irreducible $d$-tuples. In Example~\ref{ex:BKK}, the irreducible pair $(P_1,P_2)$ with the normalized mixed volume $4$ can be embedded into a maximal irreducible pair $(P_1, 2 P_1)$ with the same normalized mixed volume (see Fig. \ref{ex:BKK:enlarging}). We refer to \rs{sat} for more details. 

\begin{figure}[h]
\input{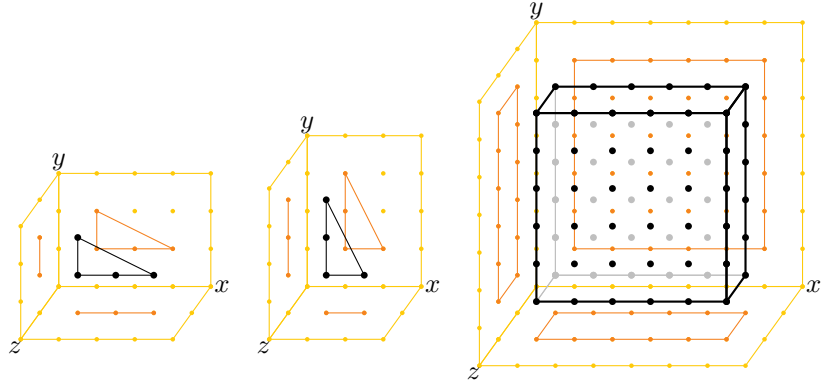}
\caption{An example of a non-irreducible triple in dimension three.\label{fig:ex:non-irredicuble}}
\end{figure}

\begin{ex}
	\label{ex:non-irred}
	Consider the triple $(P_1,P_2,P_3)$ of lattice polytopes in $\R^3$, where $P_1, P_2$ are  $2$-dimensional polytopes $P_1=\conv(0, 2 e_1, e_2)$, $P_2 = \conv(0, e_1, 2 e_2)$ in $\R^2\times\{0\}$, and $P_3$ is a $3$-dimensional lattice polytope that has width $1$ in the direction $e_3$ and is contained in the slab $\R^2 \times [0,1]$ of width $1$. The triple $(P_1,P_2,P_3)$ is not irreducible since the sum $P_1+ P_2$ is two-dimensional. The normalized mixed volume of $(P_1,P_2,P_3)$ is the product of the normalized mixed volume of the pair $(P_1,P_2)$, which is equal to $4$  (see Example~\ref{ex:BKK}), multiplied by the width $1$ of the polytope $P_3$. See also Fig.~\ref{fig:ex:non-irredicuble}
	
	This calculation can also be interpreted in light of Theorem~\ref{thm:BKK}. The triple $(P_1,P_2,P_3)$ corresponds to a generic system $f_1(x,y)=f_2(x,y)=f_3(x,y,z)=0$ in three unknowns $x,y,z$, where $f_1(x,y)=f_2(x,y)=0$ is a sub-system depending only on $x$ and $y$, with $f_1$ and $f_2$ chosen as in Example~\ref{ex:BKK}. Since the equation $f_3(x,y,z)=0$ is linear in $z$, it can written as 
	\[
		f_3(x,y,z) :=a(x,y) z + b(x,y)  = 0
	\]
	for some $a,b\in\C[x,x^{-1},y,y^{-1}]$. As explained in Example~\ref{ex:BKK}, the generic system $f_1(x,y) = f_2(x,y)=0$ has $4$ solutions in
	$(\C \setminus \{0\})^2$. If $a(x,y)$ and $b(x,y)$ are generic, then by plugging these four solutions into $f_3(x,y,z)=0$ we arrive at linear equations in $z$, each having exactly $1$ solution $z \in \C \setminus \{0\}$. Thus, the generic system $f_1(x,y) = f_2(x,y) = f_3(x,y,z) = 0$ has $4$ solutions in total.  
\end{ex}

\subsection{Our contribution} In this paper we present an algorithmic approach to Classification~Problem~\ref{class:prob} when $m$ and $d$ are given and $d \le 3$.  
Our algorithm has produced all  maximal irreducible pairs of lattice polytopes in $\R^2$ of mixed volume up to 10 and
all maximal irreducible triples of lattice polytopes in $\R^3$ of mixed volume up to $4$. 

Esterov and Gusev's approach for showing finiteness of irreducible $d$-tuples  $(P_1,\dots,P_d)$ with mixed volume $m$ is to provide
a bound on the  volume of the Minkowski sum $P_1 + \dots + P_d$  in terms of $m$, 
$$\Vol(P_1 + \dots + P_d) \leq b_d(m).$$ 
This yields a naive algorithm as, by a well-known result
of Lagarias and Ziegler \cite[Theorem~2]{LagariasZiegler}, there exists a bounding box
containing, up to lattice equivalence, all lattice polytopes of volume 
at most $b_d(m)$ in dimension $d$. Consequently, one may assume that
$P_1 + \dots + P_d$ and, hence, each $P_i$ for 
$i \in [d]$ is contained in such a bounding box and carry out an
exhaustive search for $d$-tuples that actually have normalized mixed volume $m$.
While a sharp upper bound $b_d(m)$ is not known in general, it has to be 
at least $(m+d-1)^d$, as this value is realized by the 
$d$-tuple $(\Delta_d,\dots,\Delta_d,m \Delta_d)$ of normalized mixed volume $m$, where $\Delta_d$ denotes the standard $d$-dimensional simplex.
In particular, $b_3(4)\geq 6^3=216$. Clearly, the above approach of showing the finiteness does not lead to a computationally feasible algorithm.  Already the task of enumeration of all $3$-dimensional lattice polytopes with normalized volume at most $b_3(4)$ is hopelessly hard, as the number of such polytopes is tremendous.  

In order to obtain a feasible algorithm,  we make an extensive
use of the theory of mixed volumes and mixed surface area measures. In particular we use Aleksandrov-Fenchel inequality to produce upper bounds for the normalized mixed volumes
$V(P_i,P_j,P_k)$ for all choices $i,j,k \in [3]$ and the normalized volumes of the $P_i$. Our algorithm is recursive in $m$. It has been crucial to find an appropriate case distinction that would allow to keep the enumeration procedure computationally tractable. At the highest level, we distinguish between the full-dimensional case, in which all $P_i$ are three-dimensional, and the non-full-dimensional case, in which some of the $P_i$ are two-dimensional. 

We remark that we have been able to prove that the sharp upper bound $b_d(m)$ is, in fact, $(m+d-1)^d$ for $d=2,3$ 
in the case of full-dimensional polytopes. In addition, our enumeration verifies that  $b_3(m)=(m+2)^3$ for $1\leq m\leq 4$ for
{\it all} irreducible triples. (For $d=2$ the notions of irreducible and full-dimensional coincide.) Based on this evidence we conjecture that
$b_d(m)=(m+d-1)^d$ for all irreducible $d$-tuples in an arbitrary dimension. More results on this conjecture can be found in~\cite{ABS20}.

The following result produced by our algorithm provides the answer to Classification~Problem~\ref{class:prob}
for $d=3$ and $1\leq m\leq 4$. Parts (1) and (2) of Theorem~\ref{thm:answer:to:class:problem} and Corollary~\ref{cor:answer:to:class:problem} reduce the problem to the case of $d=2$ and have already been obtained in  \cite{EsterovGusev2016}. There we use projections $\pi_{1,2} : \R^3 \to \R^2$ and $\pi_3 : \R^3 \to \R^1$ given by $\pi_{1,2}(x_1,x_2,x_3) = (x_1,x_2)$ and $\pi_{3}(x_1,x_2,x_3) = x_3$, respectively. 

\begin{thm}[Classification of triples of normalized mixed volume at most $4$]
	\label{thm:answer:to:class:problem}
	
Let  $m \in \{1,2,3,4\}$. Then
$(P_1,P_2,P_3)$ is a triple of lattice polytopes in $\R^3$ of normalized mixed volume $m$ if and only if, up to equivalence of triples, it satisfies one of the following conditions.   
	\begin{enumerate}
		\item 
For some $m_1, m_2 \in \{1,2,3,4\}$ satisfying $m=m_1 m_2$, $P_3$ is a lattice segment $\{0\}^2 \times [0,m_1]$, while 
$\pi_{1,2}(P_1)\subseteq Q_1$ and $\pi_{1,2}(P_2)\subseteq Q_2$ for some pair $(Q_1,Q_2)$ appearing in the list of 
pairs of lattice polytopes of normalized mixed volume $m_2$ given in Appendix~\ref{appendix:b}. 		
		\item 
For some $m_1, m_2 \in \{1,2,3,4\}$ satisfying $m=m_1 m_2$, $\pi_3(P_3) = [0,m_1]$, while $P_1\subseteq Q_1 \times \{0\}$, 
$P_2\subseteq Q_2\times \{0\}$ for some pair $(Q_1,Q_2)$ appearing in the list of 
pairs of lattice polytopes of normalized mixed volume $m_2$ given in Appendix~\ref{appendix:b}.
		\item 
$P_1 \subseteq Q_1, P_2 \subseteq Q_2$, and $P_3 \subseteq Q_3$ for some triple $(Q_1,Q_2,Q_3)$ appearing in the list of 
triples of lattice polytopes of normalized mixed volume $m$ given in Appendix~\ref{appendix:a}.
	\end{enumerate}  
\end{thm}

Recall that a monomial change of variables for Laurent polynomials is given by 
$(x_1,\dots, x_d)\mapsto (x^{u_1},\dots, x^{u_d})$, where $x^{u_i}=x_1^{u_{1i}}\cdots x_d^{u_{di}}$ for some unimodular matrix
$U=(u_{ij})\in\GL(d,\Z)$. We call two systems $f_1=\dots=f_d=0$ and $f_1'=\dots=f_d'=0$ {\it monomially equivalent} if after a possible permutation of the $f_i$ there is a monomial change of variables which transforms $f_i$ to $x^{a_i}f'_i$ for some monomial 
$x^{a_i}$ for every $1\leq i\leq d$. Note that monomially equivalent systems have the same number of solutions in $(\C \setminus \{0\})^d$. As an immediate consequence of Theorem~\ref{thm:answer:to:class:problem} and Theorem~\ref{thm:BKK} we obtain the following result about generic Laurent polynomial systems in three variables.

\begin{cor}[Classification of trivariate Laurent polynomial systems with at most $4$ solutions]
	\label{cor:answer:to:class:problem}
	Let $m \in \{1,2,3,4\}$. Then $f_1=f_2=f_3=0$ is a generic system of Laurent polynomials with support $(A_1,A_2,A_3)\subset(\Z^3)^3$ and which has $m$ solutions in $(\C \setminus \{0\})^3$ if and only if, up to monomial equivalence,  it satisfies one of the following conditions.   
	\begin{enumerate}
		\item 
For some $m_1, m_2 \in \{1,2,3,4\}$ satisfying $m=m_1 m_2$, $A_3\subseteq\{0\}^2 \times\{0,\dots,m_1\}$, while 
$\pi_{1,2}(A_1)\subseteq Q_1$ and $\pi_{1,2}(A_2)\subseteq Q_2$ for some pair $(Q_1,Q_2)$ appearing in the list of 
pairs of lattice polytopes of normalized mixed volume $m_2$ given in Appendix~\ref{appendix:b}. 		
		\item 
For some $m_1, m_2 \in \{1,2,3,4\}$ satisfying $m=m_1 m_2$, $\pi_3(A_3)\subseteq\{0,\dots,m_1\}$, while $A_1\subseteq Q_1 \times \{0\}$,  $A_2\subseteq Q_2\times \{0\}$ for some pair $(Q_1,Q_2)$ appearing in the list of 
pairs of lattice polytopes of normalized mixed volume $m_2$ given in Appendix~\ref{appendix:b}. 
		\item 
$A_1 \subseteq Q_1, A_2 \subseteq Q_2$ and $A_3 \subseteq Q_3$ for some triple $(Q_1,Q_2,Q_3)$ appearing in the list of 
triples of lattice polytopes of normalized mixed volume $m$ given in Appendix~\ref{appendix:a}. 
	\end{enumerate}  
\end{cor}

Examples~\ref{ex:BKK} and \ref{ex:non-irred} provide an illustration to Theorem~\ref{thm:answer:to:class:problem} 
 and Corollary~\ref{cor:answer:to:class:problem}. The triple $(P_1,P_2,P_3)$  from Example~\ref{ex:non-irred} is classified by Case~(1) of our result with $m_1=1$, $m_2=4$ and the pair $(Q_1,Q_2)$ coinciding with the pair $(\conv(0,2 e_1, e_2), 2 \conv(0, 2 e_1, e_2))$ up to equivalence.

One of the outcomes of our algorithm is the following quantitative result.

\begin{thm}
\label{thm:class_main}
Let $N_{3}(m)$ (resp. $N_{3}'(m)$) be the number of equivalence classes of triples of all (resp. $3$-dimensional) lattice polytopes in $\R^3$
of mixed volume $m$ that are irreducible and maximal. We have the following table of values for $1\leq m\leq 4$:
\smallskip
\setlength{\tabcolsep}{12pt}
\renewcommand{\arraystretch}{1.2}
\begin{center}
{\rm
\begin{tabular}{|p{.5cm}|p{1.2cm}|p{1.2cm}|}
\hline
\ $m$ & \ \ $N_{3}(m)$ &\ \ $N_{3}'(m)$\\
\hline
\ 1 & \ \ \ \ 1 & \ \ \ \ 1\\
\ 2 & \ \ \ \ 7 & \ \ \ \ 4\\
\ 3 & \ \ \ 21 & \ \ \ 10\\
\ 4  & \ \ \ 92 & \ \ \ 30\\
\hline
\end{tabular}
}
\end{center}

\end{thm}

In particular, our enumeration verifies the list of maximal
irreducible triples of lattice polytopes with mixed volume 2 
proposed in \cite{EsterovGusev2015}.

For triples of 3-dimensional lattice polytopes our enumeration produces the following structural result. We have four general constructions that produce maximal triples for an arbitrary value of the mixed volume $m$ (see also Propositions~\ref{P:adding:segments} and~\ref{P:pyramid}). Our enumeration shows that for $m\leq 4$ all, but three exceptional triples, are covered by these constructions. Recall that a polytope $P$ is a \emph{combinatorial pyramid} 
if $P$ has a facet containing all but one vertex of $P$.

\begin{thm}
\label{thm:templates}
Let $(P_1,P_2,P_3)$ be a maximal triple of 3-dimensional lattice polytopes in $\R^3$
with $V(P_1,P_2,P_3) \leq 4$. Then, up to equivalence of triples, either
\begin{enumerate}
	\item[(0)]  $P_1, P_2, P_3$ are all equal to a lattice polytope $P$,
\end{enumerate}
\begin{enumerate}
	\item
	\label{item:homothetic}  there exists a lattice polytope $P$ such that
	$P_1=\alpha P$, $P_2=\beta P$, $P_3=\gamma P$ for some integers $\alpha, \beta, \gamma\geq 1$, not all the same,
	\item 
	\label{item:line_segments}
	there exists a lattice polytope $P$ and a lattice segment
	$I$ in $\R^3$, as well as integers $\alpha,\beta\geq 1$ and $\gamma\geq 0$ such that 
	\begin{align*}		
	P_1 = P + \alpha I, P_2 = P+ \beta I
	\text{ and } P_3 = P + \gamma I,
	\end{align*}
	\item
	\label{item:pyramid}  
	there exists a lattice segment $I$ and a lattice polytope $P$, which is a combinatorial pyramid with base having two edges parallel to $I$, such that
	\begin{align*}
	P_1=P_2=P \text{ and }
	P_3=P+\alpha I,
	\end{align*}
	for some integer $\alpha\geq 1$,
	\item 
	\label{item:exceptions}  
	$(P_1,P_2,P_3)$ is one of the following exceptional triples given by
		\begin{enumerate}
			\item $P_1=P_2=\conv(0,2e_1,e_2,e_3), P_3 = P_1 + [0,e_1]$, 
			\item $P_1=P_2=\conv(0,3e_1,e_2,e_3), P_3 = P_1 + [0,e_1]$, 
			\item $P_1=P_2=\conv(0,e_1,e_2,e_1+e_2+2e_3), P_3 = P_1 + [0,e_1+e_3]$.
		\end{enumerate}
\end{enumerate}
The number of equivalence classes corresponding to the above five types are presented in the following table:
\smallskip
\setlength{\tabcolsep}{12pt}
\renewcommand{\arraystretch}{1.2}
\begin{center}
{\rm
\begin{tabular}{|p{.5cm}|p{1.3cm}|p{1.3cm}|p{1.3cm}|p{1.3cm}|p{1.3cm}|}
\hline
\ $m$ & \ type (0) & \ type (1) & \ type (2) & \ type (3) & \ type (4)\\
\hline
\ 1 & \ \ \ \ 1 & \ \ \ \ 0 & \ \ \ \ 0 & \ \ \ \ 0 & \ \ \ \ 0\\
\ 2 & \ \ \ \ 3 & \ \ \ \ 1 & \ \ \ \ 0 & \ \ \ \ 0 & \ \ \ \ 0\\
\ 3 & \ \ \ \ 6 & \ \ \ \ 1 & \ \ \ \ 1 & \ \ \ \ 1 & \ \ \ \ 1\\
\ 4 & \ \ \ 17 &\ \ \ \ 5 & \ \ \ \ 3 & \ \ \ \ 3 & \ \ \ \ 2\\
\hline
\end{tabular}
}
\end{center}
\end{thm}

\begin{figure}[h]
\input{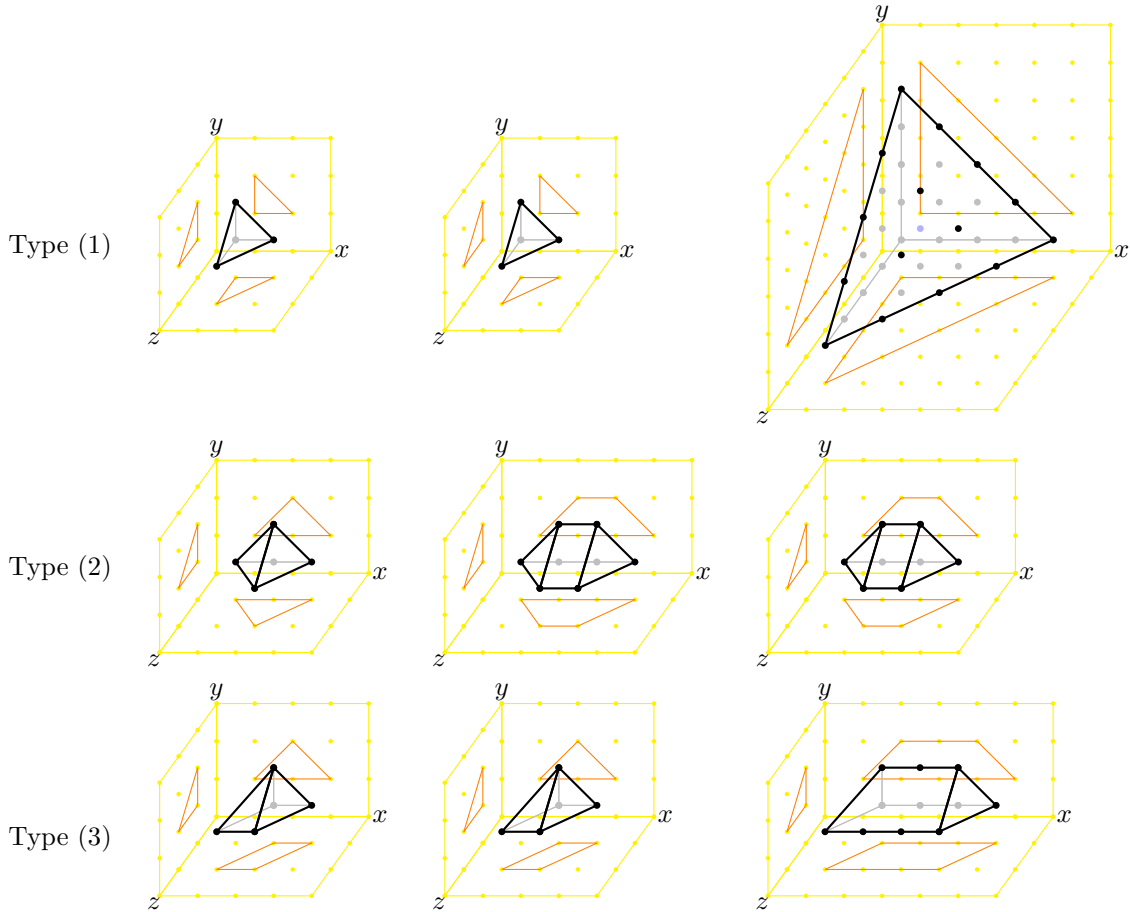}
\caption{Examples of triples of each of the types \eqref{item:homothetic}--\eqref{item:pyramid} of mixed volume 4.\label{F:templates}} 
\end{figure}

Examples of triples of polytopes of types (1)--(3) are presented in Fig.~\ref{F:templates}.
Note that we can compute the mixed volume $m$ of the triples in each of the types (0)--(3) by the following general formulas. For the triples of types (0) and (1) we have $m=\Vol(P)$ and $m=\alpha\beta\gamma\Vol(P)$, respectively. For the triples of type (2) we have
$m=\Vol(P)+ (\alpha + \beta + \gamma)\Vol_{\pi_I}(P)$, where $\Vol_{\pi_I}$ denotes the 2-dimensional
volume of the projection along the direction of $I$. Finally, for type (3) triples, we have $m=\Vol(P)+ \alpha\Vol_{\pi_I}(P)$.

Parts of our algorithm also provide a direct way to carry out enumeration of
pairs of polygons of given mixed volume. We have carried out this enumeration
for mixed volume up to 10 and obtained the following result.

\begin{thm}
Let $N_2(m)$ be the number of equivalence classes of pairs of $2$-dimensional
lattice polytopes in $\R^2$ of mixed volume $m$ that are maximal. 
We have the following table of values for $1\leq m\leq 10$.
\smallskip
\setlength{\tabcolsep}{12pt}
\renewcommand{\arraystretch}{1.2}
\begin{center}
{\rm
\begin{tabular}{|p{.5cm}|p{12mm}|}
\hline
\ $m$ & \ \  $N_2(m)$\\
\hline
\ \ 1 & \ \ \ \ \ \,1  \\
\ \ 2 & \ \ \ \ \ \,3  \\
\ \ 3 & \ \ \ \ \ \,6 \\
\ \ 4  & \ \ \ \ 13\\
\ \ 5 & \ \ \ \ 18\\
\ \ 6 & \ \ \ \ 38 \\
\ \ 7 & \ \ \ \ 46 \\
\ \ 8 & \ \ \ \ 87 \\
\ \ 9 & \ \  \,118 \\
\ 10 & \ \  \,202 \\
\hline
\end{tabular}
}
\end{center}
\end{thm}

 Our computations have been carried out using Sagemath \cite{sagemath}
 and an implementation of our enumeration algorithm
 as well as data files containing the enumeration results can
 be found at 
 \url{https://github.com/christopherborger/mixed_volume_classification}.

The approach presented in this paper also admits a natural extension
to higher dimensions. However, such an 
extension will be computationally 
much more expensive already in dimension $4$. 

 As auxiliary tools, we have also developed equivalence tests for configurations of lattice polytopes in Section~\ref{sec:equivalences} as well as
  an algorithm for enumeration of $d$-dimensional lattice polytopes with a given volume in Subsection~\ref{sect:sandwich}. The algorithmic problem of enumeration of lattice polytopes with a given volume is also interesting in its own right. Independently, another algorithm for solving this problem has recently been developed by Gabriele Balletti \cite{balletti2018volenum}.
   One of the nice features of our approach to enumeration by volume is its simplicity. Furthermore,  the general template of our volume enumeration algorithm can be modified to solve further enumeration problems in the theory of lattice polytopes that are similar in nature. All in all, we hope that the technical parts will also be useful in other contexts and will help advance the algorithmic theory of lattice polytopes, which has been emerging in the last decade. 

\subsection*{Acknowledgements.}  

This paper emerged from a discussion started at the Oberwolfach Mini-Workshop  \emph{Lattice Polytopes: Methods, Advances and Applications} (ID 1738c), which was organized by Takayuki Hibi, Akihiro Higashitani, Katharina Jochemko, and Benjamin Nill in September 2017. We would like to thank the organizers and participants of this workshop for sharing their insights and creating an open and inspiring research atmosphere. 

The first and second author were funded by the Deutsche Forschungsgemeinschaft (DFG, German Research Foundation) - 314838170, GRK 2297 MathCoRe, which also supported the third author's research visit at Otto-von-Guericke Universit\"at Magdeburg in September 2018. 

The authors would like to thank Alexander Esterov for advertising the problem and giving
valuable feedback on our enumeration results. Finally, the authors are grateful to the anonymous referees for providing a detailed list of comments and suggestions which improved the exposition. 

\section{Basic notions and background results}

\subsection{Basic notation.} Let $\N$ be the set of positive integers. For $k \in \Z_{\geq 0}$ we use $[k]$ to denote the set $\{1,\ldots,k\}$ with $[0]$ being the empty set. Throughout, $d \in \N$ denotes the dimension of the ambient space $\R^d$. We use $\subseteq$ and $\subset$ for inclusion and strict inclusion, respectively. 

The problems that we consider involve a choice of a real $d$-dimensional Euclidean space $E$ together with a lattice  $\Lambda\subset E$ of full rank. This choice is not important for our results, so most of the time we 
work with the \emph{integer lattice} $\Z^d$ in $\R^d$.
We let $e_1,\ldots,e_d$ denote the standard basis vectors in $\R^d$ which also form a basis of the lattice $\Z^d$. The notation $\Id_k$ stands for the $k \times k$ identity matrix.

If $X, Y$ are subsets of 
a vector space
$E$, then their \emph{Minkowski sum} $X+Y$ is given by $X+Y=\setcond{x+y}{x \in X, \ y \in Y}$. Furthermore, for $\lambda \in \R$ and $X \subseteq E$, we use the notation $\lambda X = \setcond{\lambda x}{x \in X}$. 

\subsection{Groups of unimodular transformations and their action.}\label{groups} 
Let $\Lambda\subset E$ be as above. We denote by $\GL(\Lambda)$ the group of all linear transformations $\phi$ on $E$ that satisfy $\phi(\Lambda) = \Lambda$. We call elements of $\GL(\Lambda)$ \emph{linear unimodular transformations}. By choosing a basis in $\Lambda$ we can identify 
$\GL(\Lambda)$  with the group of $d \times d$ \emph{unimodular matrices}, which
are the matrices $U \in \Z^{d \times d}$ with $|\det(U)| = 1$.

On several occasions, we use the group of linear unimodular transformations $\phi \in \GL(\Lambda)$ that act identically on a  $(d-1)$-dimensional linear subspace $H$ of $E$ spanned by $d-1$ linearly independent lattice vectors. We call such transformations \emph{unimodular shearings} along $H$. In particular, for $\Lambda=\Z^d$ unimodular shearings $\phi \in \GL(\Lambda)$ along the coordinate hyperplane $H=\R^{d-1} \times \{0\}$ have the form $\phi(x) = U x$ with 
\[
	U = \begin{pmatrix} \Id_{d-1} & t \\ 0 & 1 \end{pmatrix}
\] 
and $t \in \Z^{d-1}$. 

By $\Aff(\Lambda)$ we denote the group of affine transformations $\phi$ on $E$ satisfying $\phi(\Lambda) = \Lambda$. These are transformations $\phi$ of the form $\phi(x) = \psi(x) +v$, where $\psi$ is a linear unimodular transformation and $v \in \Lambda$ is a lattice vector. We call elements of $\Aff(\Lambda)$ \emph{affine unimodular transformations}. Clearly, $\GL(\Lambda)$ is a subgroup of $\Aff(\Lambda)$ and $\Aff(\Lambda)$ acts on subsets of $E$ so that $\phi \in \Aff(\Lambda)$ sends $X \subseteq E$ to its image $\phi(X)$. With this in mind, we introduce equivalences of subsets of $E$ modulo subgroups of $\Aff(\Lambda)$. If $X$ and $Y$ are subsets of $E$ and $G$ is a subgroup of $\Aff(\Lambda)$ we say that $X$ and $Y$ coincide modulo $G$ if $\phi(X) = Y$ for some $\phi \in G$. In this case, we write $X \equiv Y \mod G$.  
We also let $\Aff(\Lambda)$ act on $k$-tuples $(X_1,\ldots,X_k)$ of subsets of $E$ by $\phi( X_1,\ldots, X_k) = (\phi(X_1),\ldots,\phi(X_k))$. Thus, analogously, we say that two such $k$-tuples $(X_1,\ldots,X_k)$ and $(Y_1,\ldots,Y_k)$ coincide modulo $G$ if $(\phi(X_1),\ldots,\phi(X_k)) = (Y_1,\ldots,Y_k)$ holds for some $\phi \in G$.

There is yet another group, which we denote by $\GG_{d,k}$ 
that we let act on $k$-tuples $(X_1,\ldots,X_k)$. An element of $\GG_{d,k}$ is determined by a linear unimodular transformation $\phi \in \GL(\Lambda)$, $k$ lattice vectors $v_1,\ldots,v_k \in \Lambda$ and a permutation $\sigma$ on $[k]$. Elements of $\mathbb{G}_{d,k}$ are transformations of $(\R^d)^k$ of the form $(x_1,\ldots,x_k) \mapsto (\phi(x_{\sigma(1)}) + v_{\sigma(1)},\ldots, \phi(x_{\sigma(k)}) + v_{\sigma(k)})$. We say that two $k$-tuples $(X_1,\ldots,X_k)$ and $(Y_1,\ldots,Y_k)$ of subsets of $\R^d$ coincide up to $\GG_{d,k}$-equivalence, if $\psi(X_1,\ldots,X_k) = (Y_1,\ldots,Y_k)$ holds for some $\psi \in \GG_{d,k}$. 
We write $(X_1,\ldots,X_k) \equiv (Y_1,\ldots,Y_k) \mod \GG_{d,k}$ to denote this equivalence relation. 


\subsection{Lattice polytopes.} \label{subsect:conv} Let $E$ be a $d$-dimensional Euclidean 
space with a full-rank lattice $\Lambda$.
By $\conv$ we denote the convex hull operation on $E$. A {\it polytope} in $E$ is the convex hull of finitely many points in $E$.
For a polytope $P\subset E$ we use $\aff(P)$ to denote the affine span of $P$, that is the smallest affine subspace of $E$ containing $P$. By definition, $\dim(P)=\dim(\aff(P))$. Given a polytope $P\subset E$, its {\it support function} 
$h_P:\R^d\to\R$ is defined by 
\[
	h_P(\xi)=\max\{\langle \xi, x\rangle : x \in K\}.
\] 
Here $\langle \xi, x\rangle$ denotes the inner product in $E$ which in the case of  $E=\R^d$ is the standard inner product. 
Given $\xi\in\R^d$ we use $P^\xi$ to denote the {\it face of $P$ corresponding to $\xi$} defined
by
\[
P^\xi=\{ x \in P: \la \xi, x\ra = h_P(\xi)\}.
\]
Clearly, when $\xi=0$ we get $P^\xi=P$. For $\xi$ not equal to zero, $P^\xi$ depends only on the direction of $\xi$. 
Faces of dimension $\dim(P)-1$ are called {\it facets} and faces of dimension 0 are called {\it vertices} of $P$.  
We let $\vertset(P)$ denote the vertex set of $P$. When $\xi$ has Euclidean length $1$ and $P^{\xi}$ is a facet, $\xi$ is uniquely determined by $P^{\xi}$ and is called the {\it unit outer facet normal} of $P^\xi$. 

Given $X \subseteq E$, we denote by $\cP(X)$ the family of all non-empty polytopes $P$ with $\vertset(P) \subseteq X$. For $k \in [d]$, we use $\cP_k(X)$ to denote the family of $k$-dimensional polytopes belonging to $\cP(X)$. In this paper, we address enumeration and classification problems within the family $\cP(\Lambda)$ of \emph{lattice polytopes}. 

Let $\{v_1,\dots, v_d\}$ be a lattice basis for $\Lambda$. We call the simplex $\conv(0,v_1,\dots,v_d)$ and any of its lattice translates a
{\it unimodular simplex}. In particular, when $\Lambda=\Z^d$ and $\{e_1,\ldots,e_d\}$ is the standard basis for $\Z^d$,
we call $\Delta_d : = \conv (0,e_1,\ldots,e_d)$ the \emph{standard simplex}. Note that a simplex in $\R^d$ is unimodular
(with respect to the lattice $\Z^d$) if and only if it is $\Aff(\Z^d)$-equivalent to $\Delta_d$.

For lattice polytopes $P\in\cP(\Z^d)$
it is more convenient to consider the support function on the set $\cS^{d-1}\subset\Z^d$ of primitive vectors. 
Recall that an integer vector $u\in\Z^d$ is {\it primitive} if its entries do not have a common divisor greater than~$1$. Then
the support function $h_P(u)=\max\{\langle u, x\rangle : x \in P\}$ takes integer values for all  $u\in\cS^{d-1}$.
As before, we will speak about faces $P^u$ of $P$. In the case when $P^u$ is a facet, $u$ is called the {\it primitive outer facet normal} of $P^u$.

One of the central functionals that we want to consider on $\cP(\Lambda)$ is the Euclidean volume (i.e., the $d$-dimensional Lebesgue measure restricted to $\cP(\Lambda)$). Working with lattice polytopes, along with the {Euclidean volume} $\vol$, 
one also considers the so-called \emph{normalized volume $\Vol$ relative to $\Lambda$}, scaled so that the volume of a unimodular simplex equals one. The advantage of the scaling 
is that  $\Vol(P) \in \Z$ for each $P \in \cP(\Lambda)$. For $P\in\cP(\Z^d)$ we have $\Vol(P)=d!\vol(P)$ as
$\vol(\D_d)=1/d!$.

\subsection{Mixed volumes.} 
Let $E$ be a $d$-dimensional Euclidean space. There exists a uniquely defined functional 
\[
	v : 
\underbrace{\cP(E) \times \ldots \times \cP(E)}_d \to \R,
\] with $v(P_1,\ldots,P_d)$ being invariant under permutations of $P_1,\ldots,P_d \in \cP(E)$, such that the equality
\[
	\vol(\lambda_1 P_1 + \cdots+ \lambda_k P_k) = \sum_{i_1=1}^k \cdots \sum_{i_d=1}^k \lambda_{i_1} \cdots \lambda_{i_d} v(P_{i_1},\ldots,P_{i_d})
\]
holds for all $P_1,\dots,P_k\in\cP(E)$,
non-negative scalars $\lambda_1,\ldots, \lambda_k \ge 0$, and $k \in \N$ (see \cite[Theorem and Definition~5.1.7]{Schneider2014}). One can extend $v$ to the set of $d$-tuples of non-empty compact convex sets, but for the purposes of this paper, it will be enough to consider the case of polytopes. The value $v(P_1,\ldots,P_d)$ is called the \emph{Euclidean mixed volume} of the $d$-tuple $(P_1,\ldots,P_d)$. Replacing the Euclidean volume in the above definition with the normalized volume relative to some lattice $\Lambda\subset E$ we obtain the {\it normalized mixed volume} $V(P_1,\dots,P_d)$ relative to $\Lambda$. In what follows when we say ``volume'' or ``mixed volume'' we always assume ``normalized volume'' or  ``normalized mixed volume'' relative to some lattice. When the lattice is not specified it is meant to be 
$\Z^d\subset\R^d$ for appropriate $d\in\N$.

The mixed volume satisfies a number of properties that will be useful for our considerations. 
Their proofs can be found in
\cite[Sections 5.1, 7.3]{Schneider2014} and \cite[p. 120]{Ewald96}.

\begin{prop}\label{P:properties} For $P_1,\ldots,P_d,Q_1,\ldots,Q_d \in \cP(\R^d)$ and non-negative $\lambda,\mu\in\R$, we have
\begin{enumerate}[label=\arabic*.]
	\item $V(P_1,\ldots,P_d)\geq 0$.
	\item $V(P_1,\dots,P_d)=V(Q_1,\dots,Q_d)$, whenever $(P_1,\ldots,P_d) \equiv (Q_1,\ldots,Q_d) \mod \GG_{d,d}$.
	\item $V(\lambda P_1+\mu Q_1,P_2,\dots,P_d)=\lambda V(P_1,P_2,\dots,P_d)+\mu V(Q_1,P_2,\dots,P_d)$.
	\item $V(P_1,\ldots,P_d) \in \Z$, whenever $P_1,\ldots,P_d \in \cP(\Z^d)$.
	\item Inclusion-exclusion formula
	\begin{equation}\label{eq:MV}
V(P_1,\dots,P_d)=\frac{1}{d!}\sum_{k=1}^d(-1)^{d+k}\!\sum_{i_1<\dots<i_k}\Vol(P_{i_1}+\dots+P_{i_k}).
\end{equation}
	\item Aleksandrov--Fenchel Inequality
	$$V(P_1,P_2,P_3\ldots,P_d)^2\geq V(P_1,P_1,P_3,\ldots,P_d)V(P_2,P_2,P_3,\ldots,P_d).$$ 
	\item Monotonicity
	\begin{equation}\label{eq:monotone}
V(P_1,\dots,P_d)\leq V(Q_1,\dots,Q_d),\ \text{ whenever}\ P_1\subseteq Q_1,\dots, P_d\subseteq Q_d.
\end{equation}
	
\end{enumerate}	
\end{prop}

\begin{defn}\label{D:irr}
We say that a $d$-tuple $(P_1,\dots P_d)\in \cP(E)^d$ is {\it non-degenerate} 
if for every $I\subseteq [d]$ with $1 \le |I|\leq d$ the dimension of $\sum_{i\in I}P_i$ 
is at least $|I|$. 

We say that a $d$-tuple $(P_1,\dots P_d)\in \cP(E)^d$ is {\it irreducible} if for every $I\subseteq [d]$ with $1 \le |I|<d$ the dimension of $\sum_{i\in I}P_i$ is at least $|I|+1$. 
\end{defn}

It has been observed by Minkowski \cite[Theorem 5.1.8]{Schneider2014} that $V(P_1,\dots P_d)$ is positive if and only if the 
 $d$-tuple  $(P_1,\dots P_d)$ is non-degenerate. The notion of irreducible $d$-tuples is related to the following decomposition
 property of the mixed volume. 
 
 \begin{prop}[{\cite[Theorem~5.3.1]{Schneider2014}}] \label{P:decomp} Let $P_1,\dots, P_d \in \cP(\R^d)$ 
 be such that, for some $k \in [d]$,
 $P_1,\dots,P_k$ are contained in a rational linear subspace $L\subset\R^d$ of dimension $k$, and let $\pi_L:\R^d\to\R^d/L$
 be the projection along $L$. Then 
 $$V(P_1,\dots, P_k,P_{k+1},\dots, P_d)=V(P_1,\dots, P_k)V(\pi_L(P_{k+1}),\dots,\pi_L(P_d)).$$
 \end{prop}
 In the equality above, $V(P_1,\dots, P_k)$ is the mixed volume  of $(P_1,\ldots,P_k) \in \cP(L)^k$ relative to the sublattice
 $L\cap\Z^d$, while $V(\pi_L(P_{k+1}),\dots,\pi_L(P_d))$ is the mixed volume of  $(\pi_L(P_{k+1}),\dots,\pi_L(P_d)) \in \cP(\R^d/L)^{d-k}$ relative to the quotient lattice $\Z^d/(L\cap\Z^d)$. \rp{decomp} allows one to reduce the problem of classifying all $d$-tuples $(P_1,\dots,P_d)\in\cP(\Z^d)^d$ with given mixed volume to classifying irreducible such tuples.

The following result by Esterov paves the way for solving Classification~Problem~\ref{class:prob} via enumeration.

\begin{thm}[{\cite[Theorem~1.7]{Esterov2018}}]\label{T:finiteness}
Given $m\in\N$ there exist finitely many irreducible $d$-tuples $(P_1,\dots, P_d)\in\cP(\Z^d)^d$ with 
$V(P_1,\dots, P_d)=m$, up to $\GG_{d,d}$-equivalence.
\end{thm}

Note that this statement fails for non-irreducible tuples. For a simple example one can take $A=\conv(0,e_2)$ and
$B_n=\conv(0,e_1,ne_2)$ for $n\in\N$. Then, by \rp{decomp}, $V(A,B_n)=1$ regardless of $n$.

\subsection{Surface area measure and mixed area measure.}
 Let $P \subset \R^d$ be a polytope with unit outer facet normals $\xi_1,\dots, \xi_N$ and support function $h_P$. We have  
  \begin{equation}\label{eq:sum}
\vol(P)=\frac{1}{d}\sum_{i=1}^N h_P(\xi_i) \vol(P^{\xi_i}).
\end{equation}
Indeed, assuming $0$ lies in the interior of $P$, the above sum represents the volume of $P$ as the sum 
of the volumes of pyramids over the facets of $P$. This motivates the notion of the {\it Euclidean surface area measure} $s_P$ on 
 the sphere $\mS^{d-1}$ with finite support  $\{\xi_1,\dots,\xi_N\}\subset\mS^{d-1}$ and values $s_P(\xi_i)=\vol(P^{\xi_i})$,
  the $(d-1)$-dimensional Euclidean volume of the facet $P^{\xi_i}$, for $1\leq i\leq N$. Thus, \re{sum} can be written as
 
\begin{equation}\label{eq:vol}
\vol(P)=\frac{1}{d}\int_{\mS^{d-1}} h_P(\xi) ds_P(\xi).
\end{equation}
 We refer to \cite[Section 5.1]{Schneider2014} for details. 
 
\begin{figure}[h]
\begin{tikzpicture}[scale=0.7]
	\draw[rounded corners=.3mm, line width=1pt] (0,0) -- (2,0) -- (0,2) -- cycle;
	\node at (0.6,0.6) {\small $P$};
	\foreach \x in {-1,...,3}
	\foreach \y in {-1,...,3}
	{
		\fill (\x,\y) circle (0.08);
	}
\end{tikzpicture}
\vspace{2ex}
\\
\begin{tikzpicture}[scale=0.7]
	\draw[->,very thick,red,rounded corners=.3mm, ] (0,0) -- (0.7,0.7) node[right,black]{\small $s_P( \frac{1}{\sqrt{2}},\frac{1}{\sqrt{2}})=2 \sqrt{2}$};
\draw[->,very thick,red,rounded corners=.3mm, ] (0,0) -- (0,-1) node[below,black]{\small $s_P( 0,-1)=2$};
\draw[->,very thick,red,rounded corners=.3mm, ] (0,0) -- (-1,0) node[left,black]{\small $s_P(-1,0) = 2$};
\end{tikzpicture}
\hspace{10ex}
\begin{tikzpicture}[scale=0.7]
\draw[->,very thick,blue,rounded corners=.3mm, ] (0,0) -- (1,1) node[right,black]{\small $S_P( 1,1)=2$};
\draw[->,very thick,blue,rounded corners=.3mm, ] (0,0) -- (0,-1) node[below,black]{\small $S_P( 0,-1)=2$};
\draw[->,very thick,blue,rounded corners=.3mm, ] (0,0) -- (-1,0) node[left,black]{\small $S_P(-1,0) = 2$};
\end{tikzpicture}
\caption{Surface area measure $s_P$ and the normalized surface area measure $S_P$ of the triangle $P=\conv(0,2 e_1, 2 e_2)$\label{fig:surf:area}. While $s_P$ provides the Euclidean edge length for each outer unit edge normal vector, $S_P$ provides the lattice length for each outer primitive edge normal vector.}
\end{figure}
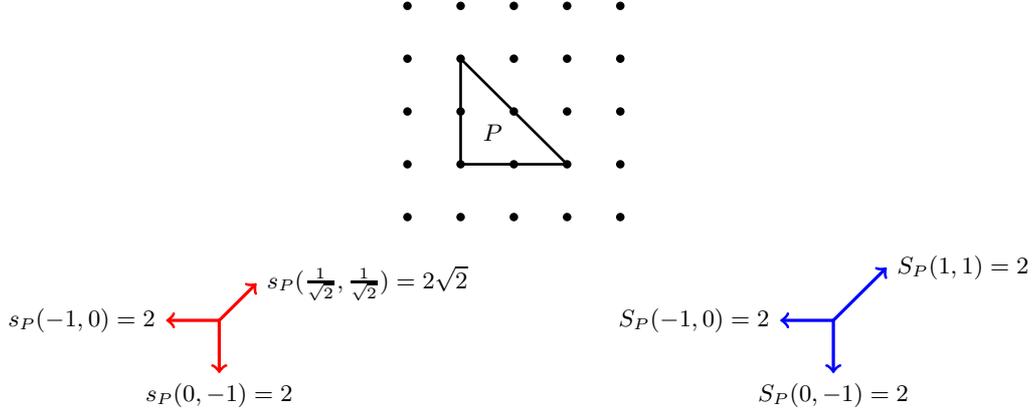

  For lattice polytopes $P\in\cP(\Z^d)$ one modifies the definition as follows. Let $\{u_1,\dots,u_N\}$ be primitive outer facet normals of $P$. 
Then we have
 \begin{equation}\label{eq:sum-Z}
\Vol(P)=\sum_{i=1}^N h_P(u_i) \Vol(P^{u_i}),
\end{equation}
where $\Vol(P^{u_i})$ is the $(d-1)$-dimensional  volume relative to the lattice $\aff(P^{u_i})\cap\Z^d$. The advantage of \re{sum-Z} is that
all the terms in the sum are non-negative integers. For this reason we introduce the {\it (normalized) surface area measure} $S_P$ on the
set of primitive vectors $\cS^{d-1}$ with support $\supp S_P=\{u_1,\dots,u_N\}$ and values $S_P(u)=\Vol(P^{u})$. Then 
\re{sum-Z} can be written as
 \begin{equation}\label{eq:integral-Z}
\Vol(P)=\int_{\cS^{d-1}} h_P(u) dS_P(u).
\end{equation}
See also Fig.~\ref{fig:surf:area} for an illustration.

Let $P_1,\dots, P_{d-1}$ be polytopes in $\R^d$. The {\it Euclidean mixed area measure} $s_{P_1,\dots,P_{d-1}}$ is defined 
by the inclusion-exclusion formula \cite[page 281]{Schneider2014}:
$$s_{P_1,\dots, P_{d-1}}=\frac{1}{(d-1)!}\sum_{k=1}^{d-1}(-1)^{d-1-k}\sum_{i_1<\dots<i_k}s_{P_{i_1}+\dots+P_{i_k}}.$$
Then, by \re{MV} in \rp{properties}, we have
\begin{equation}
	\label{eq:mixed:area:direction}
	s_{P_1,\dots, P_{d-1}}(\xi)=v(P_1^\xi,\dots, P_{d-1}^\xi)
\end{equation}
for any $\xi\in\mS^{d-1}$. This implies that $s_{P_1,\dots, P_{d-1}}$ has support
consisting of vectors $\xi\in\mS^{d-1}$ such that $(P_1^\xi,\dots, P_{d-1}^\xi)$ is non-degenerate (see the remark after \rd{irr}). 
Note that the support is contained in the set of outer facet normals of $P_1+\dots+P_{d-1}$, and hence is finite. 

Furthermore, let $\supp s_{K_1,\dots, K_{d-1}}=\{\xi_1,\dots, \xi_N\}$. We have mixed analogs of \re{sum} and \re{vol} (see \cite[Theorem 5.1.7]{Schneider2014}):

\begin{equation}\label{eq:msum}
v(P_1,\dots, P_{d-1},P)=\frac{1}{d}\sum_{i=1}^Nh_P(\xi_i)v(P_1^{\xi_i},\dots, P_{d-1}^{\xi_i}),
\end{equation}
\begin{equation}\label{eq:mvol}
v(P_1,\dots, P_{d-1},P)=\frac{1}{d}\int_{\mS^{d-1}}h_P(\xi)ds_{P_1,\dots, P_{d-1}}(\xi)
\end{equation}
for any polytope $P\subset\R^d$. 
Now let $u_i$ be the primitive vector in the direction of $\xi_i$, for $1\leq i\leq N$. Normalizing \re{msum} we obtain
\begin{equation}\label{eq:msum-Z}
V(P_1,\dots, P_{d-1},P)=\sum_{i=1}^Nh_P(u_i)V(P_1^{u_i},\dots, P_{d-1}^{u_i}).
\end{equation}
Similarly to the unmixed case we define the {\it (normalized) mixed area measure} $S_{P_1,\dots, P_{d-1}}$ on the 
set of primitive vectors $\cS^{d-1}$ with support $\supp S_{P_1,\dots, P_{d-1}}$ and values $S_{P_1,\dots, P_{d-1}}(u)=V(P_1^{u},\dots, P_{d-1}^{u})$. As before, the terms in \re{msum-Z} are non-negative integers whenever $P$ is also a lattice polytope. We write \re{msum-Z} as
 \begin{equation}\label{eq:mintegral-Z}
V(P_1,\dots, P_{d-1},P)=\int_{\cS^{d-1}}h_P(u)dS_{P_1,\dots, P_{d-1}}(u).
\end{equation}
Note that, similarly to the mixed volume, the mixed area measure is Minkowski linear in each of its arguments.

The following proposition relates the notions of irreducible tuples and mixed area measure.

\begin{prop}\label{P:irr-mixed-area-measure}
A tuple $(P_1,\dots,P_{d-1})\in\cP(\R^d)^{d-1}$ can be extended to an
 irreducible $d$-tuple if and only if $\supp s_{P_1,\ldots,P_{d-1}}$ positively spans $\R^d$.
\end{prop}

\begin{proof} Let $s$ denote the measure $s_{P_1,\ldots,P_{d-1}}$.
First note that 
\begin{equation}\label{eq:centroid}
\sum_{\xi\in \supp s}\xi s(\xi)=0.
\end{equation}
Indeed, by \re{mvol} for any $x\in\R^d$ and a polytope $P$ we have
$$0=v(P_1,\dots,P_{d-1},P+x)-v(P_1,\dots,P_{d-1},P)=\frac{1}{d}\int_{\mS^{d-1}}\langle \xi,x\rangle ds(\xi).$$
Since $x$ is arbitrary we obtain \re{centroid}. Therefore, $\supp s$ positively spans $\R^d$ if and only if
$\supp s$ linearly spans $\R^d$.

Now suppose $\supp s$ is contained in $\mS^{d-1}\cap H$ for some $(d-1)$-dimensional linear subspace $H$. 
Then, by \re{mvol}, for any segment $I$ orthogonal to $H$ we have $v(P_1,\dots,P_{d-1},I)=0$ since $h_I(\xi)=0$ for any $\xi\in \supp s$. This means that the $d$-tuple $(P_1,\dots,P_{d-1},I)$ is degenerate, which is impossible if $(P_1,\dots,P_{d-1})$ can be extended to an irreducible $d$-tuple.

Conversely, suppose $(P_1,\dots,P_{d-1})$ cannot be extended to an irreducible $d$-tuple. Without loss of generality we may assume that 
$P_1+\dots+P_k$ is contained in a $k$-dimensional subspace $L$ for some $1\leq k\leq d-1$. For
any $\xi\in \supp s$, 
$$s_{P_1,\ldots,P_{d-1}}(\xi)=v(P_1^\xi,\dots,P_{d-1}^\xi)>0,$$
and so $(P_1^\xi,\dots,P_{d-1}^\xi)$ is non-degenerate.
In particular, 
$$k\leq \dim(P_1^\xi+\dots+P_{k}^\xi)=\dim(P_1+\dots+P_{k})^\xi\leq\dim(P_1+\dots+P_k)=k.$$ 
Therefore, $P_1+\dots+P_k=(P_1+\dots+P_{k})^\xi$ which implies $\xi\in L^\perp$. Thus, $\supp s$ is a subset of
$\mS^{d-1}\cap L^\perp$.
\end{proof}

\begin{rem} One readily sees that, for tuples of lattice polytopes, the above proposition can be restated as follows: A tuple 
$(P_1,\dots,P_{d-1})\in\cP(\Z^d)^{d-1}$ can be extended to an irreducible $d$-tuple if and only if $\supp S_{P_1,\ldots,P_{d-1}}$ positively spans $\R^d$.
\end{rem}

\section{$\Z$-maximal and $\R$-maximal tuples.}\label{S:sat} 
Recall that the mixed volume $V(P_1,\ldots,P_d)$ is monotonic with respect to inclusion, see \re{monotone}. This motivates the notion of $\Z$-maximal and $\R$-maximal $d$-tuples of polytopes. 
In the definition below $R$ denotes either $\Z$ or $\R$.

\begin{defn}\label{D:sat}   Let $i \in [d]$. A $d$-tuple $(P_1,\ldots,P_d) \in \cP(R^d)^d$ is called \emph{$R$-maximal in $P_i$}, if for all
$Q_i \in \cP(R^d)$ with $P_i\subseteq Q_i$ the equality
\begin{equation}
	\label{eq:mv:strict:ineq}
	V(P_1,\ldots,P_{i-1},P_i,P_{i+1},\ldots,P_d) = V(P_1,\ldots,P_{i-1},Q_i,P_{i+1},\ldots,P_d)
\end{equation}
implies $P_i=Q_i$. We call $(P_1,\ldots,P_d)$ \emph{$R$-maximal} if it is $R$-maximal in each of the polytopes $P_1,\ldots,P_d$.  
\end{defn}

\begin{figure}[h]
	\begin{tikzpicture}[scale=0.6]
	\draw[rounded corners=.3mm, dashed,thick] (0,0) -- (4.5,0) -- (0,3) -- cycle;
	\draw[rounded corners=.3mm, line width=1pt] (0,0) -- (4,0) -- (3,1) -- (0,3) -- cycle;
	\node at (2,-2){$A=\conv(0,4 e_1, 3 e_1 + e_2, 3 e_2)$};
	\foreach \x in {-1,...,5}
	\foreach \y in {-1,...,4}
	{
		\fill (\x,\y) circle (0.08);
	}
	\end{tikzpicture}
	\hspace{10ex}
	\begin{tikzpicture}[scale=0.6]
	\draw[rounded corners=.3mm, dashed,thick] (0,0) -- (4,0) -- (0,8/3) -- cycle;
	\draw[rounded corners=.3mm, line width=1pt] (0,0) -- (4,0) -- (1,2) -- (0,2) -- cycle;
	\node at (2, -2) {$B=\conv(0, 4 e_1, e_1 + 2 e_2, 2 e_2)$};
	\foreach \x in {-1,...,5}
	\foreach \y in {-1,...,4}
	{
		\fill (\x,\y) circle (0.08);
	}
	\end{tikzpicture}
	\caption{A $\Z$-maximal pair $(A,B)$ of lattice polygons for which both $A$ and $B$ can be enlarged to non-lattice polygons $A'$ and $B'$ such that $(A',B')$ has the same mixed volume as $(A,B)$. 
	 \label{fig:Z:not:R:ex2}}
\end{figure}
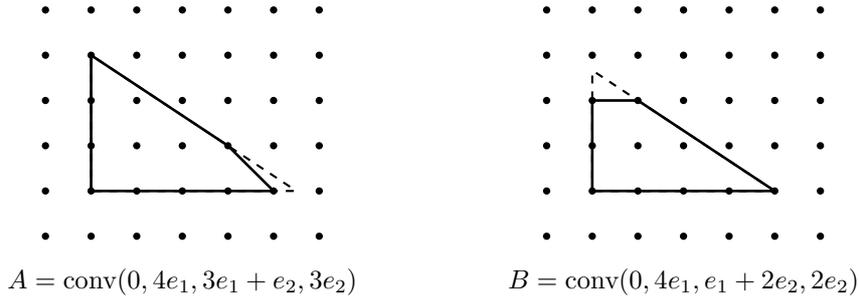

In view of the inclusion $\Z \subset \R$, if a $d$-tuple of lattice polytopes is $\R$-maximal, then it is also $\Z$-maximal. The converse is not true in general as Figs.~\ref{fig:Z:not:R:ex1} and \ref{fig:Z:not:R:ex2} illustrate.

\begin{figure}[h]
	\begin{tikzpicture}[scale=0.7]
	\draw[rounded corners=.3mm, line width=1pt] (0,0) -- (2,0) -- (0,1) -- cycle;
	\node at (1, -2) {$A=\conv(0,2 e_1, e_2)$};
	\foreach \x in {-1,...,3}
	\foreach \y in {-1,...,2}
	{
		\fill (\x,\y) circle (0.08);
	}
	\end{tikzpicture}
	\hspace{10ex}
	\begin{tikzpicture}[scale=0.7]
	\draw[rounded corners=.3mm, dashed,thick] (0,0) -- (3,0) -- (0,3/2) -- cycle;
	\draw[rounded corners=.3mm, line width=1pt] (0,0) -- (3,0) -- (1,1) -- (0,1) -- cycle;
	\node at (1.5, -2) {$B=\conv(0,3 e_1, e_1 + e_2, e_2)$};
	\foreach \x in {-1,...,4}
	\foreach \y in {-1,...,2}
	{
		\fill (\x,\y) circle (0.08);
	}
	\end{tikzpicture}
	\caption{A pair $(A,B)$ of lattice polygons, which is $\R$-maximal in $A$ and $\Z$-maximal but not $\R$-maximal in  $B$. The dashed line depict how $B$ can be enlarged to a non-lattice polygon $B' = \conv(0 , 3 e_1, \frac{3}{2} e_2)$ such that the pair $(A,B' )$ has the same mixed volume as $(A,B)$. \label{fig:Z:not:R:ex1}}
\end{figure}
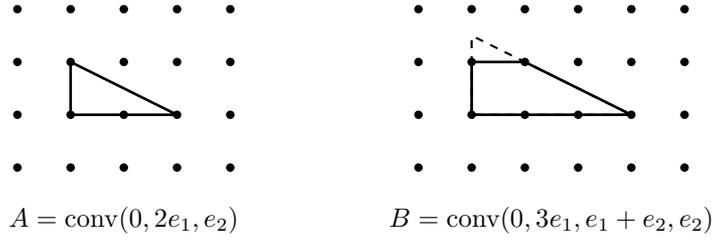

\begin{rem}\label{R:poset}
There is a natural poset structure on the set $\cP(\Z^d)^d / \GG_{d,d}$ of  $\GG_{d,d}$-equivalence classes of $d$-tuples of polytopes in $\cP(\Z^d)$. 
Let $\mathbf P$ and $\mathbf Q$ be classes containing $(P_1,\dots, P_d)$ and $(Q_1,\dots, Q_d)$, respectively. We say that
$\mathbf P\leq\mathbf Q$ if and only if $\psi(P_1,\dots, P_d)$ is contained in $(Q_1,\dots, Q_d)$ component-wise for some
$\psi \in\GG_{d,d}$. One can directly verify that this relation defines a partial order on $\cP(\Z^d)^d/ \GG_{d,d}$.
Given $m\in\N$, one may restrict this poset to the sub-poset of  $\GG_{d,d}$-equivalence classes of $d$-tuples $(P_1,\dots,P_d)$ 
of normalized mixed volume $m$. Maximal elements of this poset correspond to $\Z$-maximal $d$-tuples of polytopes in the sense of
\rd{sat}. Note that if $(P_1,\dots,P_d)$ is not irreducible then the corresponding class is not a maximal element
of the poset. Indeed, suppose $P_1,\dots, P_k$ lie in a $k$-dimensional subspace $L$, for some $1\leq k<d$. Then, by
\rp{decomp}, adding any segment whose direction vector lies in $L^\perp$ to $P_{k+1},\dots, P_d$ does not change the mixed volume of 
$(P_1,\dots,P_d)$ and, hence, the class of $(P_1,\dots,P_d)$ lies in an infinite chain in the poset. 
This shows that maximal elements are, in particular, irreducible. In what follows, however, we will always write ``maximal irreducible tuples'' for clarity. \rt{finiteness} asserts that the subposet of classes of irreducible $d$-tuples of lattice polytopes with fixed mixed volume is finite and, in particular, has finitely many maximal elements.
\end{rem}

The following proposition describes irreducible $d$-tuples $(P_1,\dots, P_{d})$ that are $\Z$-maximal in~$P_d$. In particular, it provides an algorithmic way of finding all possible $P_d\in\cP(\Z^d)$ such that $V(P_1,\dots,P_{d})=m$ and $(P_1,\dots, P_{d})$ is $\Z$-maximal in $P_d$, given a value of $m$ and a $(d-1)$-tuple $(P_1,\dots, P_{d-1})\in(\cP(\Z^d))^{d-1}$ such that $\supp s_{P_1,\ldots,P_{d-1}}$ positively spans $\R^d$ (see \ralgo{max_third}).

\begin{prop}
	\label{P:Z-maximal} Let
	$(P_1,\dots, P_{d})\in(\cP(\Z^d))^d$ be an irreducible tuple which is $\Z$-maximal in $P_d$. 
	Let $\{u_1,\dots,u_r\}$ be the support of the mixed area measure $S_{P_1,\dots,P_{d-1}}$. Then
	\begin{equation}\label{eq:1}
		P_d=\conv\{x\in \Z^d : \la u_i,x\ra\leq h_i, i\in[r]\}
	\end{equation}		
	for some $h_1,\dots,h_r\in \Z_{\geq 0}$ satisfying
	\begin{equation}\label{eq:2}
		\sum_{i=1}^rh_iS_{P_1,\dots,P_{d-1}}(u_i)=V(P_1,\dots,P_{d}).
	\end{equation}
\end{prop}
\begin{proof}
	Let $h_i=h_{P_d}(u_i)$, the value of the support function of $P_d$ at $u_i$. Then \re{2} follows directly from \re{mintegral-Z}.
	Consider 
	$$Q=\{x\in\R^d : \la u_i,x\ra\leq h_i \ \text{for all} \ i\in[r]\}.$$
	This is a rational polytope. (Its boundedness follows from \rp{irr-mixed-area-measure}.) Clearly, $P_d\subseteq Q$ and $h_{P_d}$ coincides with $h_Q$ on the support of $S_{P_1,\dots,P_{d-1}}$. Therefore, by \re{mintegral-Z} and \re{2},
	$$V(P_1,\dots, P_{d-1},Q)=\sum_{i=1}^rh_iS_{P_1,\dots,P_{d-1}}(u_i)=V(P_1,\dots,P_{d}).$$
	Now the $\Z$-maximality in $P_d$ and the monotonicity of the mixed volume imply that $P_d$ must contain all lattice points of $Q$, i.e. \re{1} holds.
\end{proof}

\begin{rem}
Note that the values for
the $h_i$ in Proposition~\ref{P:Z-maximal} are given by $h_{P_d}(u_i)$, that is the values of the support function of $P_d$ at
the support vector $u_i$ of the mixed area measure
$S_{P_1,\dots,P_{d-1}}$. This fact can be used to determine
whether a given $d$-tuple $(P_1,\dots,P_d) \in (\cP(\Z^d))^d$ is $\Z$-maximal in $P_d$ by checking whether
the equality $P_d = \conv\{x\in \Z^d : \la u_i,x\ra\leq h_{P_d}(u_i), i\in[r]\}$ holds.
\end{rem}

A similar argument provides the corresponding statement for tuples $(P_1,\dots, P_{d})\in(\cP(\R^d))^d$ that are $\R$-maximal in $P_d$.

\begin{prop}
	\label{P:R-maximal} Let
	$(P_1,\dots, P_{d})\in(\cP(\R^d))^d$ be an irreducible tuple which is $\R$-maximal in $P_d$. 
	Let $\{\xi_1,\dots,\xi_r\}$ be the support of the mixed area measure $s_{P_1,\dots,P_{d-1}}$. Then
	\begin{equation}
		P_d=\{x\in \R^d : \la \xi_i,x\ra\leq h_i \ \text{for all} \ i\in[r]\}
	\end{equation}		
	for some $h_1,\dots,h_r\in \R_{\geq 0}$ satisfying
	\begin{equation}
		\sum_{i=1}^rh_is_{P_1,\dots,P_{d-1}}(\xi_i)=v(P_1,\dots,P_{d}).
	\end{equation}
\end{prop}

For $\R$-maximality we also have the following simple criterion based on comparing the mixed area measures. Since the mixed area measures of $d$-tuples of polytopes can be computed using \eqref{eq:mixed:area:direction}, \rp{R-maximal:criterion} gives a computational test for $\R$-maximality of $d$-tuples which we utilize in our algorithms.

\begin{prop}	
	\label{P:R-maximal:criterion}
	Let $(P_1,\ldots,P_d)\in \cP(\R^d)^d$ be an irreducible tuple and assume that $P_d$ is $d$-dimensional. Then $(P_1,\ldots,P_d)$ is $\R$-maximal in $P_d$ if and only if $\supp s_{P_d} \subseteq \supp s_{P_1,\ldots,P_{d-1}}$. 
\end{prop}

 \begin{proof} To simplify notation, let $s_1=\supp s_{P_d}$ and $s_2=\supp s_{P_1,\ldots,P_{d-1}}$.
 Suppose $s_1\subseteq s_2$. Consider a polytope $Q_d\in\cP_d(\R^d)$ such that
 $P_d\subseteq Q_d$ and $$v(P_1,\dots,P_{d-1},P_d)=v(P_1,\dots, P_{d-1},Q_d).$$
 By \re{msum} we have $h_{P_d}(\xi)=h_{Q_d}(\xi)$ for every $\xi\in s_2$.  
 Then
 $$Q_d\subseteq\bigcap_{\xi\in s_2}\{x\in\R^d : \la \xi,x\ra \leq h_{P_d}(\xi)\}\subseteq\bigcap_{\xi\in s_1}\{x\in\R^d : \la \xi,x\ra \leq h_{P_d}(\xi)\}=P_d,$$ 
 since $s_1$ coincides with the set of outer facet normals of $P_d$. (Here we use the assumption that $P_d$ is $d$-dimensional.) Therefore $Q_d=P_d$.
  
  Conversely, let $(P_1,\dots, P_d)$ be $\R$-maximal in $P_d$. Consider
  $$Q_d=\bigcap_{\xi\in s_2}\{x\in\R^d : \la \xi,x\ra \leq h_{P_d}(\xi)\}.$$
  Since $(P_1,\ldots,P_d)$ is irreducible, by \rp{irr-mixed-area-measure}, $s_2$ positively spans $\R^d$, hence, $Q_d$ defines a polytope.
  Clearly, $P_d\subseteq Q_d$ and $h_{P_d}(\xi)=h_{Q_d}(\xi)$ for every $\xi\in s_2$.
  Therefore, by  \re{msum}, we have 
$$
v(P_1,\dots, P_{d-1},P_d)=v(P_1,\dots, P_{d-1},Q_d).$$
This implies that $P_d=Q_d$, and, hence, the set $s_1$ of primitive facet normals of $P_d$ is contained in~$s_2$. 
 \end{proof}

In the remainder of the section we show that the two constructions appearing in parts (\ref{item:line_segments}) and (\ref{item:pyramid}) of  
Theorem~\ref{thm:templates} produce $\R$-maximal triples. Clearly, 
the triples in part (\ref{item:homothetic}) are also $\R$-maximal, and one can use Propositions \ref{P:Z-maximal} and \ref{P:R-maximal} to see that  the triples in part (\ref{item:exceptions}) are $\Z$-maximal, but not $\R$-maximal.

\begin{prop}
	\label{P:adding:segments}
	Consider polytopes $P_1,P_2,P_3 \in \cP_3(\R^3)$ such that 
	$P_1 = P + \alpha I$, $P_2 = P + \beta I$, and $P_3 = P + \gamma I$, 
	where $P \in \cP(\R^3)$, $\dim(P) \ge 2$, $I \in \cP_1(\R^3)$, 
	and $\alpha,\beta,\gamma$ are non-negative real numbers at 
	most one of which is $0$. 
	Then the triple $(P_1,P_2,P_3)$ is $\R$-maximal. 
\end{prop}
\begin{proof}
	We verify that the triple $(P_1,P_2,P_3)$ is $\R$-maximal in 
	$P_3$ using \rp{R-maximal:criterion}. One has
	\[
		s_{P_1,P_2} = s_{P+\alpha I, P + \beta I} = s_{P} + (\alpha+\beta) s_{P,I},
	\]
	using the linearity, $s_{P,P}=s_P$, and the fact that the measure $s_I=s_{I,I}$ 	    is zero. Furthermore, one similarly obtains
	\[
		s_{P_3} = s_{P_3,P_3} = s_P + 2 \gamma s_{P,I}.
	\]
	As $\alpha,\beta,\gamma$ are nonnegative integers at most one of which
	can be $0$ this shows $\supp s_{P_3} \subseteq \supp s_{P_1,P_2}$. Showing
	$\R$-maximality in $P_1$ and $P_2$ is completely analogous.
\end{proof}

\begin{prop}\label{P:pyramid}
	Let $P' \in \cP_3(\R^3)$ be a combinatorial pyramid with base $P$ and let 
	$I$ be a segment parallel to two edges of $P$. Then the triple $(P',P',P'+ I)$ is $\R$-maximal.
\end{prop}

\begin{proof}
We are going to show that $P'$ and $P' + I$ have the same sets of facet normals
and, hence, the conditions of \rp{R-maximal:criterion} are satisfied. 
Let $P'=\conv(P\cup\{v\})$  for some $v\in\R^3$ and assume $I=[0,w]$ for some $w\in\R^3$.
Let $J_1$ and $J_2$ be two edges of $P$ parallel to $I$. 
Note that each edge of $P+I$ is equal to either an edge of $P$, or the sum of $w$ and an edge of $P$, or $J_i+ I$ for $i=1,2$. 
Similarly, each facet of $P'+I$ is equal to either a facet of $P'$, or the sum of $w$ and a facet of $P$, or
$\conv((J_i+I)\cup(v+I))$, for $i=1,2$. This implies $P'$ and $P' + I$ have the same facet normals.
\end{proof}

\section{Tests for equivalence of tuples of polytopes modulo group actions}

\label{sec:equivalences}

As we deal with enumeration up to equivalences with respect to action of various groups as described in Subsection~\ref{groups}, we need to be able to decide algorithmically whether two lattice polytopes are equivalent with respect to a given group action. More generally, we need tests for equivalence of tuples of lattice polytopes modulo group actions. In this section we present such tests for the groups of linear and affine unimodular transformations and the group $\GG_{d,k}$.

\subsection{Equivalence of polytopes modulo $\GL(\Z^d)$}

The literature contains several algorithms that test equivalence of two lattice polytopes $P,Q \in \cP_d(\Z^d)$ modulo $\GL(\Z^d)$. See, for example, \cite{KreuzerSkarke1998} and \cite{GrinisKasprzyk2013arxiv}, where the latter also provides an overview of existing techniques. The algorithm of Kreuzer and Skarke from \cite{KreuzerSkarke1998}, relying on the so-called \emph{normal form} of a lattice polytope, was implemented in {Sagemath}. The normal form  of a lattice polytope $P$ is uniquely determined by $P$. It encodes a sequence of vertices $(v_1,\ldots,v_t)$ of a polytope $\conv(v_1,\ldots,v_t)$ that coincides with $P$ up to $\GL(\Z^d)$. Two polytopes $P, Q \in \cP_d(\Z^d)$ coincide up to $\GL(\Z^d)$ if and only if their normal forms are the same. See also Example~3.4 in \cite{GrinisKasprzyk2013arxiv}. Using the normal form, each polytope in $\cP_d(\Z^d)$ can be brought into a \emph{normal $\GL(\Z^d)$-position}. In other words, in each equivalence class from $\cP_d(\Z^d) / \GL(\Z^d)$ a unique representative is chosen. Using such a normal position in enumeration algorithms is convenient because, for avoiding repetitions modulo $\GL(\Z^d)$, it suffices to bring each newly found polytope into its normal position. 

\subsection{Equivalence of polytopes modulo $\Aff(\Z^d)$.} Testing equivalence of two polytopes $P, Q \in \cP_d(\Z^d)$ modulo $\Aff(\Z^d)$ can be reduced to testing equivalence modulo $\GL(\Z^{d+1})$. Indeed,  $P, Q \in \cP_d(\Z^d)$ are equivalent modulo $\Aff(\Z^d)$ if and only if the respective pyramids $\conv(\set{0_{d+1}} \cup (P \times \{1\})), \conv(\set{0_{d+1}} \cup (Q \times \{1\})) \in \cP_{d+1}(\Z^{d+1})$ are equivalent modulo $\GL(\Z^{d+1})$. This approach works well in many cases, but for the purpose of implementation it
is more convenient to introduce a \emph{normal $\Aff(\Z^d)$-position} of polytopes in $\cP_d(\Z^d)$ by choosing a representative in each of the equivalence classes modulo $\Aff(\Z^d)$. 

Our construction is as follows.
For a polytope $P \in \cP_d(\R^d)$, consider 
\[
	c_P:= \frac{1}{|\vertset(P)|} \sum_{v \in \vertset(P)} v,
\]
which is the barycenter of the set of vertices of $P$. Furthermore, we can order points of $\R^d$ lexicographically: $x=(x_1,\ldots,x_d)$ is \emph{lexicographically smaller} than $y=(y_1,\ldots,y_d)$ if, for the smallest $i \in [d]$ with $x_i \ne y_i$, one has $x_i < y_i$. For a compact subset $X$ of $\R^d$, let $\lexmin(X)$ denote the lexicographic minimum of the set $X$.
It is not hard to see that for a polytope
$P \in \cP(\R^d)$ one has $\lexmin(P) \in \vertset(P)$. This can be shown by observing that (in notation of Section~\ref{subsect:conv}) 
$\lexmin(P)$ is obtained by first taking the face $P^{e_1}$ of $P$, then the face $(P^{e_1})^{e_2}$ of $P^{e_1}$, and so on. Continuing this iteration we eventually arrive at the unique point $( \ldots (P^{e_1})^{e_2} \ldots )^{e_d} = \lexmin(P)$.
In particular, if $P \in \cP(\Z^d)$ is a lattice polytope
then $\lexmin(P)$ is a lattice point. Based on this notions we can now introduce a normal $\Aff(\Z^d)$-position.

\begin{prop}
	\label{P:on:affine:normal:position}
	Let $P \in \cP_d(\Z^d)$ be a polytope with $N$ vertices and let $P'$ be the normal $\GL(\Z^d)$-position of the lattice polytope $N (P- c_P)$. Then the lattice polytope 
	\[
		P'' := \frac{1}{N} \left( P' - \lexmin(P') \right)
	\]
	is $\Aff(\Z^d)$-equivalent to $P$. 
\end{prop}
\begin{proof}
	Let $\phi \in \GL(\Z^d)$ be a linear unimodular transformation 
	sending $P'$ to $N(P-c_P)$. Using the fact that, for any compact
	set $X \subset \R^d$,
	one has $\lexmin(X - x) = \lexmin(X)-x$ and $\lexmin(kX)=k \lexmin(X)$ 
	for all $x \in \R^d$ and $k \in \R_{\geq 0}$ one obtains:
	\begin{align*}
	\phi(P'') 
	&= \frac{1}{N}\phi(P') - \frac{1}{N} \phi\left(\lexmin(P')\right)\\
	&= P - c_P - \phi\left(\lexmin(\phi^{-1}(P-c_P))\right) \\
	&= P - \phi\left(\lexmin(\phi^{-1}(P))\right).
	\end{align*}	 
	As $\phi(\lexmin(\phi^{-1}(P)))$ is a lattice point, this proves the claim.
\end{proof}

The polytope $P''$ in \rp{on:affine:normal:position} is uniquely determined by $P$.
We call $P''$  the \emph{normal $\Aff(\Z^d)$-position} of $P$. Note that two lattice
polytopes are equivalent modulo $\Aff(\Z^d)$ if and only if
their normal $\Aff(\Z^d)$-positions coincide.

\begin{rem}

Grinis and Kasprzyk \cite[\S3.3]{GrinisKasprzyk2013arxiv} suggest to use an affine normal form, obtained by slightly modifying the definition of a normal form of a lattice polytope.
They also observe in \cite[\S2.4]{GrinisKasprzyk2013arxiv} that $P,Q \in \cP_d(\Z^d)$ coincide modulo $\Aff(\Z^d)$ if and only if $c_P - c_Q \in \Z^d$ and $P - c_P$ and $Q - c_Q$ are equivalent modulo $\GL(\Z^d)$.  This observation is related to our \rp{on:affine:normal:position}.
\end{rem}

\subsection{Equivalence of $k$-tuples of $d$-dimensional polytopes modulo $\GG_{d,k}$}

To test equivalence of $k$-tuples $(P_1,\ldots,P_k)$ of lattice polytopes modulo $\GG_{d,k}$, we use the so-called Cayley polytope of $(P_1,\ldots,P_k)$,
which is defined as
\[
	\Cay(P_1,\dots,P_k):= \conv(P_1 \times \{e_1\} \cup \ldots \cup P_k \times \{e_k\}) \in \cP(\Z^{d}\times\Z^k).
\]
 
The following statement shows how equivalence of tuples of lattice polytopes 
modulo $\GG_{d,k}$ can be checked using the Cayley polytope construction.

\begin{prop}
	\label{P:equivalence:via:cayley}
	Let $(P_1,\ldots,P_k)$, $(Q_1,\ldots,Q_k) \in \cP(\Z^d)^k$ be $k$-tuples of non-empty 
	lattice polytopes such that $\dim(P_1 + \dots + P_k)=d$. 
	Then the following conditions are equivalent: 
	\begin{enumerate}[label={(\roman*)}]
		\item \label{item:tuple:eq}$(P_1,\dots,P_k)\equiv (Q_1,\dots,Q_k)\mod\GG_{d,k}$, 
		\item \label{item:cayley:eq}$\Cay(2 P_1,\ldots, 2 P_k)\equiv\Cay(2 Q_1, \ldots, 2 Q_k)\mod\GL(\Z^{d+k})$.
	\end{enumerate}
\end{prop}
\begin{proof}
Let $P:=\Cay(2 P_1,\ldots, 2 P_k)$ and $Q:= \Cay( 2 Q_1,\ldots, 2 Q_k)$.
The implication \ref{item:tuple:eq} $\Rightarrow$ \ref{item:cayley:eq} 
is straightforward as translating one of the polytopes $P_i$ 
corresponds to unimodular shearing the Cayley polytope along a coordinate hyperplane and applying a 
common linear unimodular transformation $\phi \colon \R^d \to \R^d$ to
each $P_i$ corresponds to applying the linear unimodular transformation
$\phi \times \Id_k \colon \R^d \times \R^k \to \R^d \times \R^k$
to the Cayley polytope. Finally, applying a permutation $\sigma$ to 
the $P_i$ corresponds to permuting the last $k$ coordinates of the
Cayley polytope. All of these transformations are in $\GL(\Z^{d+k})$. 
	
	Let us show the converse implication 
	\ref{item:cayley:eq} $\Rightarrow$  \ref{item:tuple:eq}. For a lattice polytope $F$, consider the following property $(\ast)$: for all $a,b \in \vertset(F)$, the point $(a+b)/2$ is a lattice point. This property is invariant under affine unimodular transformations. It is easy to see that a face $F$ of $P$ has property $(\ast)$ if and only if $F \subseteq (2 P_i) \times \{e_i\}$ for some $i \in [k]$. In particular, $ (2 P_i) \times \{e_i\}$ are inclusion-maximal faces of $P$ that have property $(\ast)$. Similarly, $(2 Q_i) \times \{e_i\}$ are inclusion-maximal faces of $Q$ that have property $(\ast)$. We thus see that a linear unimodular transformation $\phi$ that sends $P$ to $Q$ sends each $(2 P_i) \times \{e_i\}$, with $i \in [k]$, to some $(2 Q_{\sigma(i)}) \times \{e_{\sigma(i)}\}$, where $\sigma$ is a permutation on $[k]$. Let $\phi \colon \R^d \times \R^k \to \R^d \times \R^k$ be a linear unimodular transformation mapping $P$ to $Q$. Without loss of generality 
we may assume that $(2 P_i) \times \{e_i\}$ is mapped onto 
$(2 Q_i) \times \{e_i\}$ and that $0 \in 2P_i$
(in particular, $\phi(0,e_i)=(t_i,e_i)$ for some $t_i \in \R^d$) for all
$i \in [k]$. 
As $\phi$ maps $(2P_1 + \dots + 2P_k) \times \{e_1 + \dots + e_k\}$ onto
$(2Q_1 + \dots + 2Q_k) \times \{e_1 + \dots + e_k\}$ and
$\dim(2P_1 + \dots + 2P_k)=d$ by assumption, the map 
$\phi$ preserves the affine subspace $\R^d \times 
\{e_1 + \dots + e_k\} \subset \R^d \times \R^k$. One easily verifies
that this implies that the linear subspace $\R^d \times \set{0} \subset \R^d \times \R^k$ is also preserved. Hence, with respect to the
standard basis of $\R^{d+k}$, the map $\phi$ has the form
\begin{align*}
	\begin{pmatrix}
	U & T \\
	0 & \Id_k
	\end{pmatrix} \in \GL(\Z^{d+k})
\end{align*}
for a unimodular matrix $U \in \GL(\Z^d)$ and a matrix
$T \in \Z^{d \times k}$ whose $i$-th
column equals $t_i$ for $i \in [k]$. In particular, for any 
$i \in [k]$, the map $x \mapsto Ux + t_i$ is an affine unimodular transformation mapping $P_i$ onto $Q_i$, which proves the claim.
\end{proof}

\begin{rem}
Considering the Cayley polytopes of the second dilates
of $P_1,\dots,P_k \in \cP(\Z^d)$ is crucial. 
Consider for example the polygons $\Delta_2$ and 
$\square_2 = [0,1]^2$. Based on these polygons we construct  two pairs of $3$-dimensional 
polytopes $P_1 = P_2 =\conv(\Delta_2 \times \{0\} \cup \square_2 \times \{1\})$ and
$Q_1 =  \Delta_2 \times [0,1]$, $Q_2 =\square_2 \times [0,1]$. The pairs $(P_1,P_2)$ and $(Q_1,Q_2)$
are clearly not $\GG_{3,2}$-equivalent but 
$\Cay(P_1,P_2)$ and $\Cay(Q_1,Q_2)$ are equivalent.
\end{rem}

\section{The enumeration algorithm for full-dimensional polytopes in dimension three}

In this section we present an algorithm for solving the following
enumeration problem. Throughout this and subsequent sections the word ``maximal'' will mean ``$\Z$-maximal''.

\begin{enumprob}
	\label{mvol:full:dim}
	Given $m \in \N$, enumerate, up to $\GG_{3,3}$-equivalence, all
	 maximal triples $(P_1,P_2,P_3)\in \cP_3(\Z^3)^3$ of full-dimensional
	 lattice polytopes satisfying  $V(P_1,P_2,P_3) = m$. 
\end{enumprob}

As we mentioned in the introduction, the main idea for solving Enumeration Problem~\ref{mvol:full:dim}
is to produce upper bounds for the mixed volumes
$V(P_i,P_j,P_k)$ for all choices $i,j,k \in [3]$ 
instead of an upper bound for the volume of the Minkowski sum $P_1+P_2+P_3$. As an illustration,
consider the triple $(P_1,P_2,P_3)=(\Delta_3,\Delta_3,m \Delta_3)$ which has mixed volume $m$.
While the volumes of $P_3$ and the 
Minkowski sum $P_1+P_2+P_3$ are large, one has
$V(P_1,P_1,P_2)=1$, i.e.  some of the mixed volumes $V(P_i,P_j,P_k)$ are small. 
Such relations are enforced by the Aleksandrov-Fenchel inequality.
Proposition~\ref{P:classification_cases_3d} below characterizes 
this phenomenon in general.

\begin{prop}
\label{P:classification_cases_3d}
Let $(P_1,P_2,P_3) \in \cP_3(\Z^3)^3$ satisfy $V(P_1,P_2,P_3)=m$ for a 
given $m \in \N$. Then, up to relabeling, either $V(P_1,P_1,P_2)<m$ or $V(P_i,P_i,P_j)=m$ for all $i,j \in [3]$ with $i \neq j$. In  the latter 
case $\Vol(P_i) \leq m$ for all $i \in [3]$.
\end{prop}

\begin{proof}
Suppose there are $i,j \in [3]$ with $i \neq j$
such that $V(P_i,P_i,P_j) \neq m$. After possibly reordering we may 
assume $V(P_1,P_1,P_2) \neq m$. If $V(P_1,P_1,P_2)$ is strictly smaller 
than $m$, we have proven the claim. So let us assume $V(P_1,P_1,P_2)>m$.
In this case, by the Aleksandrov-Fenchel inequality $V(P_1,P_2,P_3)^2 \ge V(P_1,P_1,P_2) V(P_3,P_3,P_2)$ 
we have $V(P_2,P_3,P_3)<m$ so that the claim holds for the 
ordering $(P_3,P_2,P_1)$.

It is left to prove that, if $V(P_i,P_i,P_j)=m$ for all $i,j \in [3]$ with $i \neq j$,
then $\Vol(P_i) \leq m$ for all $i \in [3]$. This is a direct consequence
of the Aleksandrov-Fenchel inequality, as $V(P_i,P_i,P_j)^2 \geq V(P_i,P_j,P_j)\Vol(P_i)$ holds for any $i,j \in [3]$ with $i \neq j$.
\end{proof}

\begin{rem}\label{R:equal:polytopes}
Note that in the case $V(P_i,P_i,P_j)=m$ for all $i,j \in [3]$ with $i \neq j$ and $\Vol(P_i)=m$ for all $i \in [3]$, the 
Aleksandrov-Fenchel inequalities $V(P_i,P_i,P_j) \geq V(P_i,P_j,P_j)\Vol(P_i)$ become equalities, which implies that $P_1=P_2=P_3$.
This is due to the characterization of the equality case in Minkowski's inequality, 
see \cite[Theorem 7.2.1]{Schneider2014}.
\end{rem}

Let us now present an algorithm to solve Enumeration~Problem~\ref{mvol:full:dim}.
Note that Proposition~\ref{P:classification_cases_3d} justifies the 
restriction in Step 2 to a case a. relying on an inductive 
enumeration of maximal triples of lower mixed volume
(see Step~1) and a very specific case b.

\begin{algo}[Enumeration of full-dimensional triples]
\label{algo:full_dim_big} \
\begin{enumerate}[label={\textbf{Step~\arabic*}:}]
	\item[\textbf{Input:}] $m \in \N$.
	\item[\textbf{Output:}] A list of all
	maximal triples of full-dimensional polytopes 
	$(P_1,P_2,P_3) \in \cP_3(\Z^3)^3$ with $V(P_1,P_2,P_3)=m$,
	up to $\GG_{3,3}$-equivalence.
	\item If $m=1$, return the triple $(\Delta_3,\Delta_3,\Delta_3)$.
	Else, recursively run Algorithm~\ref{algo:full_dim_big} for
	all input values $m' < m$ in order to obtain 
	a list of all
	 maximal triples of full-dimensional polytopes $(P_1,P_2,P_3) \in \cP_3(\Z^3)^3$ with $V(P_1,P_2,P_3) < m$, up 
	to $\GG_{3,3}$-equivalence.
	\item 
	\begin{enumerate}
		\item[\textbf{a.}] Enumerate, up to $\GG_{3,2}$-equivalence,
		 all pairs $(P_1,P_2) \in \cP_3(\Z^3)^2$ with $V(P_1,P_1,P_2)<m$
		 (see Remark~\ref{rem:ind_full}).
		\item[\textbf{b.}] Enumerate, up to $\GG_{3,2}$-equivalence,  all pairs $(P_1,P_2) \in \cP_3(\Z^3)^2$ with 
		$V(P_1,P_1,P_2)=V(P_2,P_2,P_1)=m$
		(see Algorithm~\ref{algo:full_dim_sub}).
	\end{enumerate}
	\item For a given pair $(P_1,P_2) \in \cP_3(\Z^3)^2$ as in
	Step 2.a 
	enumerate, up to translations, all $P_3 \in \cP_3(\Z^3)$ such that 
	$V(P_1,P_2,P_3)=m$ and such that the triple $(P_1,P_2,P_3)$ is maximal in $P_3$ (see Algorithm~\ref{algo:max_third}). 
	\item Given a triple $(P_1,P_2,P_3) \in \cP_3(\Z^3)^3$ as in Step 2.b
	check whether it is maximal in $P_1$ and $P_2$ and, if so,
	add it modulo $\GG_{3,3}$-equivalence to the final list of
	maximal triples of mixed volume $m$. 
\end{enumerate}
\end{algo}

\begin{rem}[De-maximization procedure]
\label{rem:ind_full}
 We may obtain the
pairs of Step 2.a from the list of all maximal triples
of full-dimensional lattice polytopes $(P_1,P_2,P_3) \in \cP_3(\Z^3)^3$ with $V(P_1,P_2,P_3) < m$, as recursively obtained
in Step~1. Note that we need to 
consider not only those pairs $(P_1,P_2) \in \cP_3(\Z^3)^2$ for which
the triple $(P_1,P_1,P_2)$ is maximal, but all pairs such that 
$V(P_1,P_1,P_2)=m'<m$. These can be obtained by iteratively pealing
off vertices of the maximal triples of mixed volume less than $m$
and searching among them for triples of the form $(P_1,P_1,P_2)$ up to permutations
and translations.
The running time of this task is very reasonable for values
$m' \in \set{1,2,3}$ but is growing very fast in $m'$ and would also
be growing extensively if we were to consider higher dimensions.
\end{rem}

Dealing with Step 2.b is more involved and, hence, we employ the following
algorithm:

\begin{algo}[Step 2.b of Algorithm~\ref{algo:full_dim_big}] \
\label{algo:full_dim_sub}
\begin{enumerate}[label={\textbf{Step \arabic*}:}]
	\item[\textbf{Input:}] $m \in \N$.
	\item[\textbf{Output:}] A list of all 
	 pairs $(P_1,P_2) \in \cP_3(\Z^3)^2$ with 
	$V(P_1,P_1,P_2) = V(P_1,P_2,P_2) = m$, up to $\GG_{3,2}$-equivalence.
	\item Enumerate, up to $\Aff(\Z^3)$-equivalence, all $P_1 \in \cP_3(\Z^3)$ with $\Vol(P_1) \leq m$ (Enumeration~Problem~\ref{enum:vol}).
	\item Given $P_1 \in \cP_3(\Z^3)$ with $\Vol(P_1) \leq m$, determine, 
	up to translations, all $Q \in \cP_3(\Z^3)$ such that $V(P_1,P_1,Q)=m$ 
	and the triple $(P_1,P_1,Q)$
	is maximal in $Q$ (see Algorithm~\ref{algo:max_third}).
	\item Given a pair $(P_1,Q) \in \cP_3(\Z^3)^2$ as in Step 2, determine all subpolytopes
	$P_2 \subseteq Q$ such that $\Vol(P_2) \leq m$ and $V(P_1,P_1,P_2)=V(P_2,P_2,P_1)=m$ (see Section~\ref{S:sandwich_subpolytopes}).
	\item Given a pair $(P_1,P_2) \in \cP_3(\Z^3)^2$ with 
	$P_2 \subseteq Q$ as above, add it modulo $\GG_{3,2}$-equivalence 
	to the final list.
\end{enumerate}
\end{algo}

\section{Extension to the general case in dimension three}

In this section we present an extension of Algorithm~\ref{algo:full_dim_big}
 allowing us to enumerate general maximal irreducible triples 
 $(P_1,P_2,P_3)\in P(\Z^3)^3$
without assuming that all $P_i$ are full-dimensional, thus solving
the following enumeration problem.

\begin{enumprob}
	\label{mvol:maximal}
	Given $m \in \N$, enumerate, up to $\GG_{3,3}$-equivalence, all
	maximal irreducible triples $(P_1,P_2,P_3)\in \cP(\Z^3)^3$ satisfying
	$V(P_1,P_2,P_3) = m$.
\end{enumprob}

While we may still formulate
a statement analogous to \rp{classification_cases_3d}, the existence
of $i,j \in [3]$ with $i \neq j$ such that $V(P_i,P_i,P_j)<m$ does not necessarily
allow us to build upon the enumeration for lower mixed volumes as in
Remark~\ref{rem:ind_full}. The reason is that, if one has $\dim(P_i)=2$, 
the triple $(P_i,P_i,P_j)$ is not 
irreducible anymore and therefore may not be contained in one of the maximal
irreducible triples of lower mixed volume. Hence, we carry out a different
case distinction for triples that contain at least one
polytope of dimension $2$. For a 2-dimensional lattice polytope
$P \in \cP(\Z^3)$, let $\Vol_2(P)$ denote the normalized 2-dimensional volume relative to the affine lattice $\aff(P)\cap\Z^3$. 

\begin{prop}\label{prop:class_cases_2d}
Let $(P_1,P_2,P_3) \in \cP(\Z^3)^3$ be an irreducible triple such that $V(P_1,P_2,P_3)=m$ and at least one of the $P_i$ is $2$-dimensional. Then 
there exist distinct indices $i,j \in [3]$ such that one of the following
holds:
\begin{enumerate}
\item[\rm{(a)}] $(P_i,P_i,P_j)$ is irreducible and satisfies $V(P_i,P_i,P_j)<m$,
\item[\rm{(b)}] $\dim(P_i)=2$, $\dim(P_j)=3$, $V(P_i,P_i,P_j) \leq m$, and
$V(P_j,P_j,P_i)=m$,
\item[\rm{(c)}] $\dim(P_i)=\dim(P_j)=2$, $V(P_i,P_i,P_j) \leq m$, and 
$V(P_j,P_j,P_i) \leq m^2$.
\end{enumerate}
\end{prop}

\begin{proof}
Without loss of generality we may assume $\dim(P_1)\leq \dim(P_2) \leq \dim(P_3)$.
We distinguish cases based on the dimensions of the polytopes in the tuple.
Assume first $\dim(P_1)=2$ and $\dim(P_2)=\dim(P_3)=3$. Consider the 
Aleksandrov-Fenchel inequality $m^2=V(P_1,P_2,P_3)^2\geq V(P_2,P_2,P_1)
V(P_3,P_3,P_1)$. If $V(P_2,P_2,P_1)<m$ or $V(P_3,P_3,P_1)<m$, 
one has (a) for $(i,j)=(2,1)$ or $(i,j)=(3,1)$, respectively. Otherwise, 
one has $V(P_2,P_2,P_1)=V(P_3,P_3,P_1)=m$. Now, if (a) does not hold for $(i,j)=(2,3)$ then
$V(P_2,P_2,P_3)\geq m$ and the Aleksandrov-Fenchel inequality $m^2\geq V(P_1,P_1,P_3)V(P_2,P_2,P_3)$
implies $V(P_1,P_1,P_3)\leq m$, i.e. (b) holds for $(i,j)=(1,3)$. Similarly, we show that either (a) holds for $(i,j)=(3,2)$
or (b) holds for $(i,j)=(1,2)$. 

Assume now that $\dim(P_1)=\dim(P_2)=2$, and $\dim(P_3)=3$. 
Consider the Aleksandrov-Fenchel
inequality $m^2 \geq V(P_2,P_2,P_1) V(P_3,P_3,P_1)$. If $V(P_3,P_3,P_1)<m$,
then (a) holds for $(i,j)=(3,1)$. Otherwise, one has $V(P_2,P_2,P_1)\leq m$.
As, additionally, $V(P_1,P_1,P_2) V(P_3,P_3,P_2) \leq m^2$, in this
case (c) holds for $(i,j)=(2,1)$.

Let us finally assume that $\dim(P_1)=\dim(P_2)=\dim(P_3)=2$. Then the inequality 
$m^2 \geq V(P_1,P_1,P_2) V(P_3,P_3,P_2)$ yields that either 
$V(P_1,P_1,P_2) \leq m$ or $V(P_3,P_3,P_2) \leq m$.
Analogously to the above, one also has $V(P_2,P_2,P_1) \leq m^2$
and $V(P_2,P_2,P_3) \leq m^2$. Therefore, case (c) holds for 
$(i,j)=(1,2)$ or $(i,j)=(3,2)$.
\end{proof}

Let us now present an extension of Algorithm~\ref{algo:full_dim_big} that
allows us to solve Enumeration~Problem~\ref{mvol:maximal}. In particular,
Algorithm~\ref{algo:lower_dim_big} is used to enumerate maximal irreducible triples of a given mixed volume
containing at least one $2$-dimensional lattice polytope. Note that
Proposition~\ref{prop:class_cases_2d} justifies the restriction to the
three cases a., b., and c. in Step 2.

\begin{algo}[Extension of Algorithm~\ref{algo:full_dim_big} to general maximal irreducible triples] \
\label{algo:lower_dim_big}
\begin{enumerate}[label={\textbf{Step~\arabic*}:}]
	\item[\textbf{Input:}] $m \in \N$.
	\item[\textbf{Output:}] A list of all
	 maximal irreducible triples
	$(P_1,P_2,P_3) \in \cP(\Z^3)^3$ with $V(P_1,P_2,P_3) = m$ and
	$\dim(P_i)=2$ for at least one $i \in [3]$, 
	up to $\GG_{3,3}$-equivalence.
	\item If $m=1$, return an empty list (as the only maximal 
	irreducible triple of mixed volume $1$ is 
	$(\Delta_3,\Delta_3,\Delta_3)$ and therefore full-dimensional).
	Else, recursively run Algorithm~\ref{algo:full_dim_big} and
	Algorithm~\ref{algo:lower_dim_big} for input values
	$m' < m$ in order to obtain a list of all
	maximal irreducible triples 
	$(P_1,P_2,P_3) \in \cP(\Z^3)^3$ with $V(P_1,P_2,P_3) < m$, 
	up to $\GG_{3,3}$-equivalence.
	\item 
	\begin{enumerate}
		\item[\textbf{a.}] Enumerate, up to $\GG_{3,2}$-equivalence, 
		all pairs $(P_1,P_2) \in \cP(\Z^3)^2$ such that the triple
		$(P_1,P_1,P_2)$ is irreducible with $V(P_1,P_1,P_2)<m$
		(see Remark~\ref{rem:ind_2d}).
		\item[\textbf{b.}] Enumerate, up to $\GG_{3,2}$-equivalence,
		all pairs
		 $(P_1,P_2) \in \cP(\Z^3)^2$ with $\dim(P_1)=2$ and $\dim(P_2)=3$
		and where $V(P_1,P_1,P_2) \leq m$ and $V(P_2,P_2,P_1)=m$
		(see Algorithm~\ref{algo:lower_dim_sub}).
		\item[\textbf{c.}] Enumerate, up to $\GG_{3,2}$-equivalence,
		all pairs
		 $(P_1,P_2) \in \cP(\Z^3)^2$ where $\dim(P_1)=\dim(P_2)=2$,
		$V(P_1,P_1,P_2) \leq m$ and $V(P_1,P_2,P_2) \leq m^2$
		(see Algorithm~\ref{algo:lower_dim_sub}).
	\end{enumerate}
	\item For a given pair $(P_1,P_2) \in \cP(\Z^3)^2$ as in Step 2, 
	enumerate, up to translations, all $P_3 \in \cP(\Z^3)$ such that $V(P_1,P_2,P_3)=m$ and the triple
	$(P_1,P_2,P_3)$ is irreducible and maximal in $P_3$ (see Algorithm~\ref{algo:max_third}). 
	\item Given an irreducible triple $(P_1,P_2,P_3) \in \cP(\Z^3)^3$ as
	in Step 3, check whether it is also maximal in $P_1$ and $P_2$ and, 
	if so, add it modulo $\GG_{3,3}$-equivalence to the final list of
	maximal triples of mixed volume $m$.
\end{enumerate}
\end{algo}

\begin{rem}
\label{rem:ind_2d}
The enumeration in Step 2.a can be obtained
from the list of maximal triples of mixed volume $m'<m$
of Step~1 analogously
to the procedure described in Remark~\ref{rem:ind_full}.
\end{rem}

In order to treat Step~2.b and Step~2.c. of 
Algorithm~\ref{algo:lower_dim_big} we apply the following
algorithm. While we treat both cases in a similar way, the 
separation in some of the steps has an important computational advantage. 
This is because it allows us to have relatively small bounding boxes in
which one has to perform a rather expensive search for full-dimensional
subpolytopes (Step~2.b), while we may restrict the search inside 
 larger bounding boxes to lower-dimensional subpolytopes (Step~2.c).

\begin{algo}[Step~2.b and Step~2.c of Algorithm~\ref{algo:lower_dim_big}] \
\label{algo:lower_dim_sub}
\begin{enumerate}[label={\textbf{Step~\arabic*}:}]
	\item[\textbf{Input:}] $m \in \N$.
	\item[\textbf{Output:}] \
	\begin{enumerate}
		\item[\textbf{for b:}] A list of all pairs
		 $(P_1,P_2) \in \cP(\Z^3)^2$, up to $\GG_{3,2}$-equivalence,
		  with $\dim(P_1)=2$ and $\dim(P_2)=3$
		and where $V(P_1,P_1,P_2) \leq m$ and $V(P_2,P_2,P_1)=m$.
		\item[\textbf{for c:}]  A list of all pairs
		 $(P_1,P_2) \in \cP(\Z^3)^2$, up to $\GG_{3,2}$-equivalence,
		  with $\dim(P_1)=\dim(P_2)=2$
		and where $V(P_1,P_1,P_2) \leq m$ and $V(P_2,P_2,P_1) \leq m^2$.
	\end{enumerate}
	\item[\textbf{Step~1}:] Enumerate all $P_1 \coloneqq P'_1 \times \set{0} \in \cP(\Z^3)$
	for $P'_1 \in \cP_2(\Z^2)$ with $\Vol_2(P_1) \leq m$ up to
	equivalence (Enumeration~Problem~\ref{enum:vol}).
	\item[\textbf{Step~2}:] 
	\text{ }
	\begin{enumerate}
		\item[\textbf{for b}:] Given $P_1$ as above determine bounding boxes
	$Q_1,\dots,Q_r \in \cP_3(\Z^3)$ containing, up to shearing along the affine span of $P_1$
	and  translations, all $P_2$ which satisfy 
	$V(P_1,P_1,P_2) \leq m$ and $V(P_1,P_2,P_2)=m$ (see 
	Lemma~\ref{lem:boundingbox:lowerdim}).
		\item[\textbf{for c}:] Given $P_1$ as above determine bounding
	boxes $R_1,\dots,R_s \in \cP_3(\Z^3)$ containing, up to shearing along the affine span of $P_1$
	and  translations, all $P_2$ satisfying 
	$V(P_1,P_1,P_2) \leq m$ and $V(P_1,P_2,P_2) \leq m^2$ (see
	Lemma~\ref{lem:boundingbox:lowerdim}).
	\end{enumerate}
	\item[\textbf{Step~3}:] 
	\text{ }
	\begin{enumerate}
		\item[\textbf{for b}:] Determine all subpolytopes 
	$P_2 \in \cP_3(\Z^3)$ of the 
	bounding boxes $Q_1,\dots,Q_r$ that satisfy 
	$V(P_1,P_1,P_2) \leq m$ and $V(P_1,P_2,P_2)=m$
	(see Remark~\ref{rem:contains_segment} and Algorithm~\ref{algo:subpolytopes}).
		\item[\textbf{for c}:] Determine all subpolytopes 
	$P_2 \in \cP_2(\Z^3)$ 
	of the bounding boxes
	$R_1,\dots,R_s$ that satisfy $V(P_1,P_1,P_2) \leq m$ and 
	$V(P_1,P_2,P_2) \leq m^2$ (see Remark~\ref{rem:contains_segment} and Algorithm~\ref{algo:subpolytopes}).
	\end{enumerate}
	\item[\textbf{Step~4}:] Given $P_1$ and $P_2$ as above
	add the pair $(P_1,P_2) \in \cP(\Z^3)^2$ modulo 
	$\GG_{3,2}$-equivalence to the final list.
\end{enumerate}
\end{algo}

\section{Details of the enumeration algorithms}

In this section we provide further details of the enumeration algorithms presented in the previous sections. 

\subsection{Finding maximal $P_3$}
\label{sect:findingmaxP3}

A problem that we have to solve in various steps of the enumeration algorithm is
the following.

\begin{enumprob}
\label{enum:maximal_third} Let  $m \in \N$ and let $(P_1,P_2) \in \cP(\Z^3)^2$  be
a pair of lattice polytopes  satisfying $\dim(P_1),\dim(P_2)\geq 2$ and $\dim(P_1+P_2)=3$.
Enumerate, up to translations, all lattice polytopes
$P_3 \in \cP(\Z^3)$ such that $V(P_1,P_2,P_3)=m$ and such that the triple $(P_1,P_2,P_3)$
is irreducible and $\Z$-maximal in $P_3$.
\end{enumprob}

We solve this enumeration problem by making use of \rp{Z-maximal}.

\begin{algo}[Finding maximal $P_3$] \
\label{algo:max_third}
\begin{enumerate}[label={\textbf{Step~\arabic*}:}]
	\item[\textbf{Input:}] A pair $(P_1,P_2) \in \cP(\Z^3)^2$ with
	$\dim(P_1),\dim(P_2)\geq 2$ and $\dim(P_1+P_2)=3$ and a number $m \in \N$.
	\item[\textbf{Output:}] A list of all lattice polytopes $P_3 \in \cP(\Z^3)$
	up to translations
	such that $V(P_1,P_2,P_3)=m$ and the triple $(P_1,P_2,P_3)$ is 
	$\Z$-maximal in $P_3$.
	\item Compute the mixed area measure of $P_1$ and $P_2$.
	That is, compute the normalized mixed areas $V(P_1^{u},P_2^{u})$
	for all $u \in \Z^3$ that are primitive outer
	normal vectors of a facet of the
	Minkowski sum $P_1 + P_2$. Obtain a vector 
	$a = (a_1,\dots,a_r) \in \N^r$ of the mixed areas of $P_1,P_2$ 
	with respect to those primitive normal vectors 
	$u_1,\dots,u_r \in \Z^3$ that yield a positive mixed area.
	\item Determine all vectors $h = (h_1,\dots,h_r) \in (\Z_{\geq 0})^r$
	satisfying $\sum_{i=1}^r h_ia_i = m$.
	\item Given a vector $h \in (\Z_{\geq 0})^r$ as above, 
	compute $P_3\coloneqq\conv\set{x \in \Z^3 \colon \langle u_i , x \rangle \leq h_i\text{ for all } i \in [r]}$ and
	check whether the triple $(P_1,P_2,P_3)$ is irreducible and satisfies
	$V(P_1,P_2,P_3)=m$. If it does, append $P_3$ modulo translations to 
	the final list.
\end{enumerate}
\end{algo}

\begin{rem}
	Algorithm~\ref{algo:max_third} allows us to profit from the 
	restriction to maximal triples (or triples that are
	maximal in at least one polytope) in our enumeration. 
	For example, fixing the pair $(\Delta_3,\Delta_3) \in \cP_3(\Z^3)^2$ 
	and mixed volume $m=4$, Algorithm~\ref{algo:max_third} directly
	determines $Q = 4 \Delta_3$ as the unique maximal lattice polytope
	such that $V(\Delta_3,\Delta_3,Q)=4$.
\end{rem}

\begin{rem}
	A slight modification of Algorithm~\ref{algo:max_third} can be used in order to 
	enumerate maximal pairs of polygons $(P_1,P_2) \in \cP_2(\Z^2)$ of a given mixed 
	volume $m$.	By the Aleksandrov-Fenchel inequality in the two-dimensional case one 
	may assume without loss of generality that $\Vol(P_1)\leq m$. For any fixed 
	$P_1 \in \cP_2(\Z^2)$ one may compute the area measure and,
	analogously to Step~2 and Step~3 of Algorithm~\ref{algo:max_third}, determine a list of
	all $P_2 \in \cP_2(\Z^2)$ such that $(P_1,P_2)$ is $\Z$-maximal in $P_2$ and
	$V(P_1,P_2)=m$.
\end{rem}

\subsection{Bounding $P_2$ given a lower-dimensional $P_1$}

In the lemma below $A-A$ denotes the \emph{difference set} $A+(-A)$ of a convex set $A\subset \R^d$ and 
$A^*=\{y\in\R^d : \langle x,y\rangle\leq 1\text{ for all } x\in A\}$ denotes the \emph{polar dual} convex set. 

\begin{lem}
	\label{lem:boundingbox:lowerdim}
	Let $(P_1,P_2) \in \cP(\Z^3)^2$ be a pair of lattice polytopes such that $P_1$ is $2$-dimensional of the form   
	$P_1 = P' \times \set{0} \subset \R^2 \times \set{0}$ and $P_2$ has dimension at least 2 and positive width $w$ in the direction of $e_3$. Assume $V(P_1,P_1,P_2)=m_1$ and
	$V(P_2,P_2,P_1) \leq m_2$ for some $m_1,m_2\in\N$.
	Then, up to a shearing along $\R^2 \times \{0\}$ and a lattice translation, $P_2$ is contained in the bounding polytope
	
	\[
	R_{q_1,q_2} :=\conv \left( \setcond{\begin{pmatrix} x_1 \\ x_2 \\ x_3 \end{pmatrix} \in \Z^2 \times \{0,\ldots,w-1\}}{\begin{pmatrix}
		x_1 \\
		x_2
		\end{pmatrix}
		\in
		Q' + \frac{1}{w} 
		\begin{pmatrix}
		q_1 x_3 \\
		q_2 x_3
		\end{pmatrix} } \right),
	\]
	where $q_1, q_2 \in \{0,\ldots,w-1\}$ and 
	\[
	Q' :=\begin{pmatrix}
	0 & \frac{m_2}{w} \\
	-\frac{m_2}{w} & 0
	\end{pmatrix}
	(P'-P')^{*}.
	\]	
\end{lem}

\begin{proof}
	We may assume $0 \in P_2$ and $h_{P_2}(-e_3)=0$ and
	therefore Proposition~\ref{P:decomp} yields
	 $m_1 = V(P_1,P_1,P_2) = \Vol_2(P') h_{P_2}(e_3)$.
	Then $P_2$ contains a lattice point $(q_1,q_2,w)$ at height $w$
	and, up to shearing, we may assume that $q_1,q_2 \in \set{0,\dots,w-1}$. 
	Let $x = (x_1,x_2,x_3) \in P_2 \cap \Z^3$
	be another lattice point of $P_2$ and consider the triangle
	$T_x \coloneqq \conv((0,0,0),(q_1,q_2,w),(x_1,x_2,x_3)) \subseteq
	P_2$. Let
	$(n_1,n_2,n_3) = (q_1,q_2,w) \times (x_1,x_2,x_3)$ be the normal 
	vector to $\aff(T_x)$ with lattice length equal to $\Vol_2(T_x)$.
	Then \eqref{eq:msum-Z} yields 
	\begin{align*}
	V(T_x,T_x,P_1) &= h_{P_1}((n_1,n_2,n_3)) + h_{P_1}(-(n_1,n_2,n_3))\\
	&= h_{P'}((n_1,n_2)) + h_{P'}(-(n_1,n_2)) \\
	&= h_{P'-P'}((n_1,n_2)).
	\end{align*}
	Explicit computation of the cross product yields
	$(n_1,n_2)= (-wx_2,wx_1) + (q_2 x_3, -q_1 x_3)$.
	By the monotonicity of
	the mixed volume we have $m_2 \geq V(P_2,P_2,P_1) \geq V(T_x,T_x,P_1)$
	and, hence, 
	\begin{align*}
	(-wx_2,wx_1) + (q_2 x_3, -q_1 x_3)= (n_1,n_2) \in m_2 (P'-P')^{*}.
	\end{align*}
	This is equivalent to
	\begin{align*}
	\begin{pmatrix}
	x_1 \\
	x_2
	\end{pmatrix}
	\in
	\begin{pmatrix}
	0 & \frac{m_2}{w} \\
	-\frac{m_2}{w} & 0
	\end{pmatrix}
	(P'-P')^{*} + \frac{1}{w} 
	\begin{pmatrix}
	q_1 x_3 \\
	q_2 x_3
	\end{pmatrix},
	\end{align*}
	which shows the assertion.
\end{proof}

\begin{rem}
\label{rem:contains_segment}
Note that, in the setting of Lemma~\ref{lem:boundingbox:lowerdim}, the
bounding box $R_{q_1,q_2}$ is actually constructed under the assumption
that $P_2$ contains the segment 
$I_{q_1,q_2} = \conv((0,0,0),(q_1,q_2,w))$.
Therefore in order to enumerate the set of all suitable $P_2$ we may
restrict to searching for all $q_1,q_2 \in \set{0,\dots,w-1}$ for 
lattice polytopes inside $R_{q_1,q_2}$ that contain $I_{q_1,q_2}$.
We use this fact when we apply Algorithm~\ref{algo:subpolytopes}.
Also note that any lattice polytope $P \in \cP(\Z^3)$ with $I_{q_1,q_2} \subset P 
\subseteq R_{q_1,q_2}$ satisfies $V(P_1,P_1,P) = m_1$ by 
construction of $R_{q_1,q_2}$, 
while the upper bound of $m_2$ on $V(P,P,P_1)$ may in general
not be satisfied for some subpolytope $P \subseteq R_{q_1,q_2}$.
\end{rem}

\subsection{Enumeration of lattice polytopes by volume using sandwich factory algorithm}
\label{sect:sandwich}

\begin{enumprob}
	\label{enum:vol}
	Given $m \in \N$ and $d \in \N$, enumerate up to affine unimodular transformations all polytopes $P \in \cP_d(\Z^d)$ with $1 \le \Vol(P) \le m$.  
\end{enumprob}


\subsubsection{Sandwich-factory based approach}

We present a relatively simple algorithm to Enumeration~Problem~\ref{enum:vol} which we also found
to lead to very reasonable running times. 
The running time of the Sagemath \cite{sagemath} implementation of our algorithm on a regular desktop computer was just a few minutes for $d=3$ and $m=4$. For $d=3$, $m=6$ our implementation terminates within an hour. For $d=2$, much larger values of $m$ can be handled with an hour time limit. Even more important in the context of our original enumeration problem about mixed volumes is the fact that we use our algorithm for solving Enumeration~Problem~\ref{enum:vol} as a template for solving further similar enumeration problems by appropriately modifying the basic steps of the algorithm.   

We also refer to \cite{balletti2018volenum} for an alternative approach to enumeration of lattice polytopes by their volume. Note also that \cite{EsterovGusev2016} provides an explicit description of lattice polytopes of arbitrary dimension $d$ with the normalized volume at most $4$.  

	We call a pair $(A,B)$ of polytopes $A,B \in \cP_d(\Z^d)$ a \emph{sandwich} if $A$ is a subset of $B$. The basic principle of our algorithm is to capture all possible polytopes in a set of sandwiches $(A,B)$. If for $P \in \cP_d(\Z^d)$ the inclusion $A \subseteq P \subseteq B$ holds, we say that {\it $P$ occurs in the sandwich $(A,B)$}. The algorithm maintains a \emph{sandwich factory}, which is a set of sandwiches with the property that each $P$ in question occurs in some of the sandwiches from the set. We call the difference $\Vol(B)-\Vol(A)$ the \emph{volume gap} of a sandwich $(A,B)$. The algorithm starts with a sandwich factory containing sandwiches with a large volume gap. In the main part, the algorithm iteratively replaces sandwiches with a large volume gap by sandwiches with a smaller volume gap. Eventually, only sandwiches with volume gap $0$ remain; such sandwiches correspond to polytopes $P$ with $\Vol(P) \le m$. Thus, as soon as there are no sandwiches with positive volume gap, the enumeration task is completed. 

	\subsubsection{Initialization of the sandwich factory.}\label{S:initialization}
	It is clear that every lattice polytope $P \in \cP_d(\Z^d)$ contains an empty lattice simplex $A$, that is a simplex with exactly $d+1$ lattice points. Also, if $\Vol(P) \le m$ then, clearly, $\Vol(A) \le m$. Thus, we start with a set of sandwiches $(A,B)$ which involves all possible empty simplices $A$ with $1 \le \Vol(A) \le m$. In dimension $d=2$ there is only one empty simplex up to $\Aff(\Z^2)$-equivalence, namely, the triangle $A=\Delta_2$. In dimension $d=3$, by White's classification (see \cite{White1964} or \cite[Theorem 5]{Reznick06}), every empty $3$-dimensional simplex is $\Aff(\Z^3)$-equivalent to either the standard simplex $\Delta_3$ or $\conv(0, e_1, e_3, e_3 + a e_1 + b e_2)$ with $a,b \in \N$, $a<b$, and $\gcd(a,b)=1$. Thus, it suffices to determine such simplices $A$ with $1 \le \Vol(A) \le m$. To complete the initialization of the sandwich factory, one needs to choose an appropriate $B$ for each $A$ so that $(A,B)$ is a sandwich, which contains all lattice polytopes $P$ with $1 \le \Vol(P) \le m$ and the property $A \subseteq P$. It is intuitively clear that if a point $x$ is far away from $A$, then the volume $\conv(A \cup \{x\})$ must be large. This informal idea is expressed explicitly in the following lemma. 

\begin{lem}
	\label{lem:container}
		Let $A$ be a $d$-dimensional simplex and let $m \ge \Vol(A)$. Then 
	\[
	\setcond{x \in \R^d}{ \Vol(\conv(A \cup \{x\})) \le m} \subseteq \lambda A + (1-\lambda) c_A,
	\] 
	where $c_A$ is the barycenter of $A$ and $\lambda = (d+1) \left( \frac{m}{\Vol(A)} -1 \right) + 1$.
\end{lem}
\begin{proof}
	The proof for $d=3$ can be found in \cite[Lemma~13]{AKW17}. The proof extends directly to the case of an arbitrary dimension $d \in \N$. 
\end{proof}

In view of Lemma~\ref{lem:container}, one can fix $B$ to be the integral hull of $\lambda A + (1-\lambda) c_A$, that is
$$B=\conv ( (\lambda A + (1-\lambda) c_A ) \cap \Z^d).$$  It may still be the case that $B$ is chosen to be too large in the sense that $B$ may contain vertices $v$ with the property that $\Vol(\conv(A \cup \{v\})) > m$. Clearly, if a polytope $P \in \cP_d(\Z^d)$ occurs in $(A,B)$ 
and has the property $\Vol(P)\leq m$
then $P$ cannot contain $v$ as above. This means that one can iteratively make $B$ smaller by removing vertices $v$ as above, as long as such vertices exist. More precisely, while $v$ as above exists, one iteratively replaces $B$ by $\conv((B \cap \Z^d) \setminus \{v\})$. We call this procedure the \emph{reduction} of $B$ relative to $A$. Having carried out the above reduction of $B$ for each $A$, we complete the initialization of the sandwich factory. 

\subsubsection{Iterative updates of sandwich factory}
\begin{figure}[h]
	\begin{tikzpicture}[scale=0.4,baseline=1.5mm]
	\draw [rounded corners=.3mm, line width=1pt] (0,0) -- (1,0) -- (0,1) -- cycle;
	\draw [rounded corners=.3mm, line width=1pt] (-1,0) -- (0,-1) -- (2,-1) -- (2,0) -- (0,2) -- (-1,2) -- cycle;
	\filldraw[fill=white] (-1,2) circle (0.16);
	\foreach \x in {-2,...,3}
	\foreach \y in {-2,...,3}
	{
		\fill (\x,\y) circle (0.08);
	}
	\end{tikzpicture}
	\hspace{1em}
	$\xrightarrow{\hspace{3em}}$
	\hspace{1em}
	\begin{tikzpicture}[scale=0.4,baseline=1.5mm]
\draw [rounded corners=.3mm, line width=1pt] (0,0) -- (1,0) -- (-1,2) -- cycle;
\draw [rounded corners=.3mm, line width=1pt] (-1,0) -- (0,-1) -- (2,-1) -- (2,0) -- (0,2) -- (-1,2) -- cycle;
\filldraw[fill=white] (-1,2) circle (0.16);
\foreach \x in {-2,...,3}
\foreach \y in {-2,...,3}
{
	\fill (\x,\y) circle (0.08);
}
\end{tikzpicture}
\hspace{4em}
\begin{tikzpicture}[scale=0.4,baseline=1.5mm]
\draw [rounded corners=.3mm, line width=1pt] (0,0) -- (1,0) -- (0,1) -- cycle;
\draw [rounded corners=.3mm, line width=1pt] (-1,0) -- (0,-1) -- (2,-1) -- (2,0) -- (0,2) -- (-1,1)  -- cycle;
\filldraw[fill=white] (-1,2) circle (0.16);
\foreach \x in {-2,...,3}
\foreach \y in {-2,...,3}
{
	\fill (\x,\y) circle (0.08);
}
\end{tikzpicture}
	\caption{Replacing a sandwich by two other sandwiches with a smaller volume gap.\label{fig:sandwich:splitting}}
\end{figure}
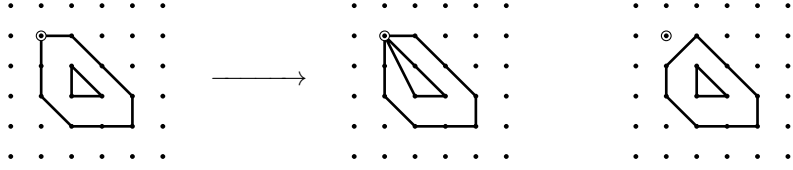

The purpose of the iterative procedure is to reduce the maximum volume gap occurring in a sandwich factory. That is, as long as there are sandwiches with positive volume gap, one considers those sandwiches $(A,B)$ in the sandwich factory, for which the volume gap $\Vol(B) - \Vol(A)$ is maximized. For each such sandwich $(A,B)$ one picks a vertex $v$ of $B$ not belonging to $A$. Every polytope $P \in \cP_d(\Z^d)$ with $\Vol(P) \le m$ occurring in $(A,B)$  may or may not contain $v$. If $P$ contains $v$, we can enclose $P$ into the sandwich $(\conv(A \cup \{v\}),B)$ with a smaller volume gap. If $P$ does not contain $v$, we can enclose $P$ into the sandwich $(A,\conv((B\cap \Z^d) \setminus \{v\}))$, whose volume gap is also smaller. Thus, we can remove the sandwich $(A,B)$ from the factory and replace it by two other sandwiches (see also Fig.~\ref{fig:sandwich:splitting}). 

Here it should also be noticed that, when we let $A$ grow, by considering $(\conv(A \cup \{v\}, B)$, we can make $B$ smaller. Indeed, $B$ may contain vertices $w$ with the property that $\Vol(\conv(A \cup \{v,w\})) > m$. Then rather than adding the sandwich $(\conv(A \cup \{v\}, B)$, we first reduce $B$ relative $\conv(A \cup \{v\})$ to a potentially smaller polytope $B'$ and then add $(\conv(A \cup \{v\}), B')$ to the sandwich factory. 

\subsubsection{Equivalent sandwiches.}
\begin{figure}[h]
	\begin{tikzpicture}[scale=0.7,baseline=-5pt]
	\draw [rounded corners=.3mm, line width=1pt] (0,0) -- (1,0) -- (0,1) -- cycle;
	\draw [rounded corners=.3mm, line width=1pt] (-1,0) -- (0,-1) -- (2,-1) -- (2,0) -- (0,2) -- (-1,2) -- cycle;
	\foreach \x in {-2,...,3}
	\foreach \y in {-2,...,3}
	{
		\fill (\x,\y) circle (0.08);
	}
	\end{tikzpicture}
	\hspace{10ex}
	\input{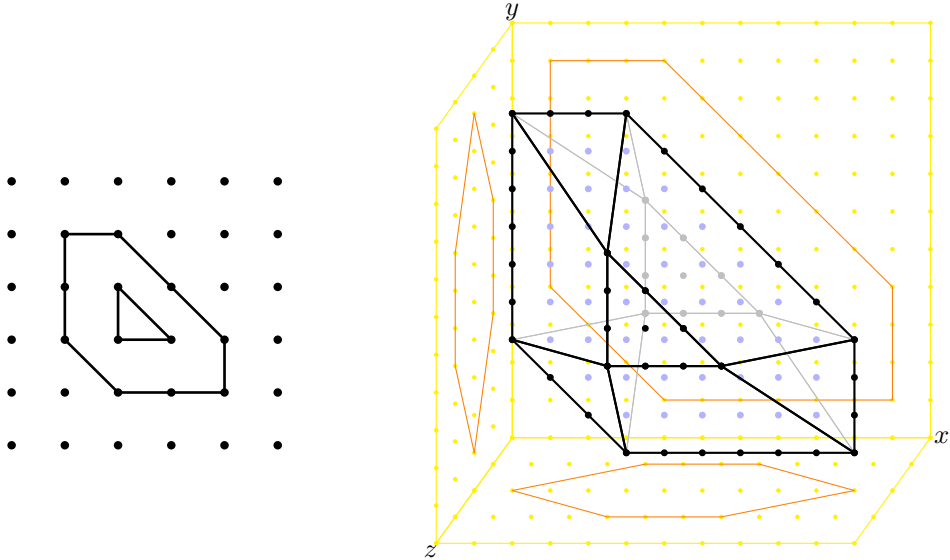}
	\caption{An example of a sandwich $(A,B)$ in dimension two (left) and the three-dimensional lattice polytope $P_{A,B}$ assigned to this sandwich (right), whose affine normal form is used to distinguish sandwiches up to affine unimodular transformations.\label{fig:layered:polytope}}
\end{figure}

While the above algorithmic steps can already be used for finding all polytopes $P \in \cP_d(\Z^d)$ with $\Vol(P) \le m$, its efficiency would not be very good as one would generate many polytopes that are equivalent up to affine unimodular transformation. When for two sandwiches $(A,B)$ and $(A',B')$ there exists an affine unimodular transformation $\phi$ with $\phi(A)=A'$ and $\phi(B) = B'$, then, up to affine unimodular transformations, the lattice polytopes occurring in $(A,B)$ also occur in $(A',B')$ and vice versa. We call such sandwiches $(A,B)$ and $(A',B')$ {\it equivalent.} Thus, if a sandwich $(A,B)$ is already present in the sandwich factory, there is no need to add $(A',B')$. Based on this idea, we add a new sandwich $(A,B)$ to the sandwich factory only if the factory does not already contain a sandwich equivalent to $(A,B)$. The test for equivalence of sandwiches can be reduced to the test for equivalence of lattice polytopes as follows. If $(A,B)$ is a sandwich then, by embedding $3 A$ into $\R^{d}\times\R$ at heights $1$ and $-1$ and $3B$ at height $0$, we obtain the polytope
\[
	P_{A,B} := \conv \bigl( \underbrace{(3 A ) \times \{1\}}_{\text{height} \ 1} \cup \underbrace{(3 B) \times \{0\}}_{\text{height} \ 0} \cup \underbrace{(3 A) \times \{-1\}}_{\text{height} \ -1} \bigr) \in \cP_{d+1}(\Z^{d}\times\Z),
\]
see also Fig.~\ref{fig:layered:polytope}.

\begin{lem}
	Two sandwiches $(A,B)$ and $(A',B')$ are equivalent if and only if the polytopes $P_{A,B}$ and $P_{A',B'}$ are equivalent up to affine unimodular transformations.
\end{lem}
\begin{proof}
	The first implication is direct. If $(A,B)$ and $(A',B')$ are equivalent,
	then there exists an affine unimodular transformation 
	$\phi \in \Aff(\Z^d)$ such that	$\phi(A)=A'$ and $\phi(B)=B'$.
	The map $\phi \times \Id \in \Aff(\Z^d \times \Z)$ 
	then satisfies $(\phi \times \Id)(P_{A,B}) = P_{A',B'}$.
	In order to show the reverse implication assume that $P_{A,B}$ and 
	$P_{A',B'}$ are equivalent and let $\psi \in \Aff(\Z^d \times \Z)$ be 
	an affine unimodular transformation such that 
	$\psi(P_{A,B})=P_{A',B'}$. Note that for both  $P=P_{A,B}$
	and $P=P_{A',B'}$ the vector $v = e_{d+1}$ is the
	unique direction such that  
	$|h_P(v)-h_P(-v)| = 2$. Therefore $\psi$ maps 
	the intersection of $P_{A,B}$ with any of the hyperplanes 
	$\R^d \times \set{-1}$, $\R^d \times \set{0}$ and 
	$\R^d \times \set{1}$ to the intersection of $P_{A',B'}$
	with the respective hyperplane. Here we use that, as $P_{A,B}$
	is symmetric with respect to the hyperplane $\R^d \times \set{0}$,
	we may assume that the intersections with the hyperplanes
	$\R^d \times \set{-1}$ and $\R^d \times \set{1}$ are not 
	permuted by $\psi$. In particular, 
	$\psi(3B \times \set{0}) = 3B' \times \set{0}$
	and, as $B$ is full-dimensional, 
	 $\psi(\R^d \times \set{h}) = \R^d \times \set{h}$ for
	any $h \in \R$. Furthermore, we may assume that $\psi$ is
	linear and, hence,  with respect to the standard basis 
	of $\R^{d+1}$ to be of the form
	\begin{align*}
	\begin{pmatrix}
	U & t \\
	0 & 1
	\end{pmatrix} \in \GL(\Z^{d+1}),
	\end{align*}
	for a unimodular matrix $U \in \GL(\Z^d)$ and
	an integer vector $t \in \Z^d$. 
	Denote by $\phi \in \GL(\Z^d)$ the linear unimodular
	transformation corresponding to $U$. Then
	$\psi((3A \times \set{-1})) = (\phi(3A)-t) \times \set{-1}$
	and 
	$\psi((3A \times \set{1})) = (\phi(3A)+t) \times \set{1}$.
	In particular $\phi(3A)-t = 3A' = \phi(3A)+t$ 
	and therefore one has $t=0$. So $\phi$ is a unimodular
	transformation such	that $\phi(A) = A'$ and 
	$\phi(B) = B'$ and hence the sandwiches $(A,B)$ and 
	$(A',B')$ are equivalent.	
\end{proof}

\subsubsection{Summary of the sandwich-factory algorithm.}
We give a complete description of the algorithm we have developed above.  

\begin{algo}[Sandwich-factory algorithm] \
\label{algo:sandwichfactory}
\begin{enumerate}[label={\textbf{Step~\arabic*:}}]
	\item[\textbf{Input:}] Dimension $d \in \N$ and volume bound $m \in \N$.
	\item[\textbf{Output:}] A list of all full-dimensional lattice polytopes
	$P \in \cP_d(\Z^d)$ with $\Vol(P) \leq m$, up to affine unimodular 
	transformations.
	\item Enumerate, up to affine unimodular transformations, all empty lattice simplices $A$ with $\Vol(A) \le m$. 
	\begin{itemize}
		\item For each $A$ as above, choose $B$ to be the integral hull 
		\[
		\conv ( (\lambda A + (1-\lambda) c_A ) \cap \Z^d),
		\] where $c_A$ is the barycenter of $A$ and 
		\[
		\lambda = (d+1) \left( \frac{m}{\Vol(A)} - 1 \right) + 1 
		\]
		and then reduce $B$ relative to $A$. 
		\item Initialize the sandwich factory $\cF$ with all pairs $(A,B)$ obtained as above.
	\end{itemize}
	\item While $\cF$ contains sandwiches with a positive volume gap, carry out the following steps for sandwiches $(A,B)$ whose volume gap is maximized: 
	\begin{itemize}
		\item pick a vertex $v$ of $B$, not contained in $A$,
		\item fix $A' := \conv(A \cup \{v\})$,
		\item compute the reduction $B'$ of $B$ relative to $A'$,
		\item fix $B'' := \conv((B \cap \Z^d) \setminus \{v\})$,
		\item add $(A',B')$ to $\cF$, unless $\cF$ already contains a sandwich equivalent to $(A',B')$,
		\item add $(A,B'')$ to $\cF$, unless $\cF$ already contains a sandwich equivalent to $(A,B'')$. 
	\end{itemize}
	\item In this step, all sandwiches $(A,B)$ in $\cF$ have the form $A=B$. Return the set of all $A$, with $(A,B) \in \cF$. Up to affine unimodular transformations, this is the set of all polytopes $P \in \cP_d(\Z^d)$ with $1 \le \Vol(P) \le m$. 
\end{enumerate}
\end{algo}

\begin{rem}
	In \rs{initialization} we described an efficient implementation of Step~1 for dimension two and three (the dimensions we are interested in, in the context of this paper). For higher dimensions we do not specify how to implement Step~1 and only notice that it can be implemented algorithmically. We also note that rather than using empty lattice simplices of normalized volume at most $m$, one can start with {\it all} lattice simplices of normalized volume at most $m$. Such simplices can be enumerated using the Hermitian normal form (see \cite[Section~4.1]{schrijver1998book}) rather easily. 
\end{rem}

\subsection{Sandwich type search for subpolytopes}
\label{S:sandwich_subpolytopes}

In this section we describe our approach towards the task of
finding all subpolytopes fulfilling certain conditions inside 
a given bounding polytope. For our purposes we found it
computationally efficient to employ an algorithm similar to
the one presented in Section~\ref{sect:sandwich}.
Our enumeration employs three different 
variations of this task that we solve using three different
variations a., b., and c. of Algorithm~\ref{algo:subpolytopes}. 
In particular, variation~a. is employed for the search for 
full-dimensional subpolytopes inside of a maximal polytope as 
obtained using Algorithm~\ref{algo:max_third}, while 
variation b. and c. deal with the search for full-dimensional or
$2$-dimensional subpolytopes of a bounding polytope as in
Lemma~\ref{lem:boundingbox:lowerdim}. Note that a sandwich type
search seems particularly natural for the search inside bounding polytopes
as in Lemma~\ref{lem:boundingbox:lowerdim}, as by 
Remark~\ref{rem:contains_segment} it suffices to search for those 
subpolytopes that contain a given segment $I$ depending on the
bounding polytope.

For a sandwich $(A,B) \in \cP(\Z^3)^2$ we call the nonnegative number
$|B \cap \Z^3| - |A \cap \Z^3|$ the \emph{lattice point gap} of $(A,B)$. Furthermore, for a lattice 
polytope $P \in \cP(\Z^3)$ let $\Vol_r(P)$ denote the normalized $k$-dimensional volume relative to the affine lattice $\aff(P)\cap\Z^3$ where $k \in [3]$ is the dimension of $P$.
Generalizing the concept of the reduction of a lattice polytope 
$B \supseteq A$ relative to $A$ used in Algorithm~\ref{algo:sandwichfactory}, we define the {\it reduction
of $B$ relative to $A$  with respect to the conditions
$\Vol_r(\,\cdot\,)\leq m_1$ and $V(\,\cdot\,,\,\cdot\,,P_1) \leq m_2$} to be 
 the polytope 
\begin{align*}
B'=\conv\{ x \in B \cap \Z^3 : A_x \coloneqq \conv(A \cup \set{x})\text{ satisfies } \Vol_r(A_x) \leq m_1,V(A_x,A_x,P_1) \leq m_2 \}.
\end{align*}
\phantom{a}
\begin{algo}[Sandwich approach to subpolytopes] \
\label{algo:subpolytopes}
\begin{enumerate}[label={\textbf{Step~\arabic*:}}]
	\item[\textbf{Input:}] \
		\begin{itemize}
			\item[\textbf{for a}:] A bounding box $M \in \cP_3(\Z^3)$,
	a lattice polytope $P_1 \in \cP(\Z^3)$, and a bound 
	$m \in \N$.
			\item[\textbf{for b/c}:] A bounding box $M \in \cP_3(\Z^3)$,
			a segment $I \subset M$, a lattice polytope
			$P_1 \in \cP(\Z^3)$, and bounds $m_1,m_2 \in \N$.
		\end{itemize}
	\item[\textbf{Output:}] \
		\begin{itemize}
			\item[\textbf{for a}:] 
			A list of all full-dimensional
			lattice polytopes $P_2$, up to translations, with 
			$P_2 \subseteq M$ such that $\Vol(P_2) \leq m$ and 
			$V(P_2,P_2,P_1) \leq m$.
			\item[\textbf{for b}:] A list of all
			 full-dimensional lattice polytopes
			$P_2$, up to translations,
			 with $I \subset P_2 \subseteq M$ such 
			that $\Vol(P_2) \leq m_1$ and $V(P_2,P_2,P_1) \leq m_2$.
			\item[\textbf{for c}:] A list of all
			$2$-dimensional lattice polytopes $P_2$, up to
			 translations, with 
			$I \subset P_2 \subset M$ such that $\Vol_2(P_2) \leq m_1$ and 
			$V(P_2,P_2,P_1) \leq m_2$.
		\end{itemize}
	\item \
		\begin{itemize}
			\item[\textbf{for a}:] Initialize the sandwich factory $\cF$
			with pairs $(S,M)$ where $S$ ranges over all empty 
			simplices in $M$ satisfying the bounding conditions
			(in particular, $V(S,S,P_1) \leq m$). Set $m_1=m_2=m$.
			\item[\textbf{for b/c}:] Initialize the sandwich factory $\cF$
			with the pair $(I,M')$, where $M'$ is the reduction of $M$ relative to $I$ with respect to the conditions
			$\Vol_r(\,\cdot\,)\leq m_1$ and $V(\,\cdot\,,\,\cdot\,,P_1) \leq m_2$. 
		\end{itemize}		 
	\item While $\cF$ contains sandwiches with positive lattice point gap, 
	carry out the following steps for sandwiches $(A,B)$ 
	having the maximal lattice point gap among the sandwiches in $\cF$: 
	\begin{itemize}
		\item pick a vertex $v$ of $B$, not contained in $A$,
		\item fix $A' := \conv(A \cup \{v\})$ (note that, as $B$ is 
		reduced relative to $A$
		with respect to the conditions $\Vol_r(\,\cdot\,)\leq m_1$ and
		$V(\,\cdot\,,\,\cdot\,,P_1) \leq m_2$, the polytope $A'$ is ensured
		to satisfy the 
		bounding conditions),
		\item (for c:) if $\dim(A')=2$ and $\dim(B)=3$, set 
		$B \coloneqq B \cap \aff(A')$,
		\item compute the reduction $B'$ of $B$ relative to $A$
		with respect to the conditions $\Vol_r(\,\cdot\,)\leq m_1$ and
		$V(\,\cdot\,,\,\cdot\,,P) \leq m_2$,
		\item fix $B'' := \conv((B \cap \Z^d) \setminus \{v\})$,
		\item add $(A',B')$ to $\cF$, unless $\cF$ already contains a translate of $(A',B')$,
		\item add $(A,B'')$ to $\cF$, unless $\cF$ already contains a translate of $(A,B'')$,
		\item remove $(A,B)$ from $\cF$. 
	\end{itemize}
	\item In this step all sandwiches $(A,B)$ have lattice point gap 0 and therefore 
	fulfill $A=B$. Return a list of $A$ for all	sandwiches $(A,B) \in \cF$.
\end{enumerate}
\end{algo}

\begin{rem}
While the overall structure of Algorithm~\ref{algo:subpolytopes}
above is very similar to Algorithm~\ref{algo:sandwichfactory} there are 
some differences. In Algorithm~\ref{algo:subpolytopes} we also
work with sandwiches $(A,B)$ for which $\dim(A) < \dim(B)$
and therefore the volume gap is not necessarily strictly decreasing in our
iterative steps. We deal with this by considering the lattice point gap
of a sandwich instead. 
Furthermore, we only identify two sandwiches $(A,B)$ and $(A',B')$
if there is a translation vector $t \in \Z^3$ such that 
$(A',B') = (A+t,B+t)$. Also note that, in addition to a volume bound
for $P_2$, we also have a bound for the mixed volume $V(P_2,P_2,P_1)$
and therefore perform a slightly different reduction step.
\end{rem}

\bibliographystyle{alpha}
\bibliography{lit}

\clearpage

\appendix

\newcolumntype{C}{>{\centering\arraybackslash} m }
\newcommand{\mycaption}[1]{\begin{center}\textbf{#1}\end{center}}
\newgeometry{left=1cm,right=1cm,bottom=1cm}

\section{Enumeration data for dimension 3}\label{appendix:a}

\begin{adjustwidth}{2cm}{2cm}
We present the complete list of maximal irreducible triples
of lattice polytopes with mixed volume at most 4. 
We omit explicitly writing down triples of the form $(P,P,P)$ for a 
lattice polytope $P \in \cP_3(\Z^3)$, that is, triples of type (0) of
Theorem~\ref{thm:templates}. Instead, we present lists of full-dimensional
lattice polytopes $P \in \cP_3(\Z^3)$ of normalized volume up to 4.
The lists of full-dimensional $\R$-maximal triples are further subdivided 
into the types (1)--(3) as in Theorem~\ref{thm:templates}.

The layout of the figures of polytopes is explained in Fig.~\ref{fig:expl:illustration}.

\begin{figure}[h]
	\begin{tikzpicture}
	\draw[->,very thick] (0.000000000000000, 0.000000000000000) -- (1.00000000000000, 0.000000000000000) node[right]{$e_1$};
	\draw[->,very thick] (0.000000000000000, 0.000000000000000) -- (0.000000000000000, 1.00000000000000) node[above]{$e_2$};
	\draw[->,very thick] (0.000000000000000, 0.000000000000000) -- (-0.500000000000000, -0.700000000000000) node[below]{$e_3$};
	\end{tikzpicture}	
	\hspace{5mm}
	\input{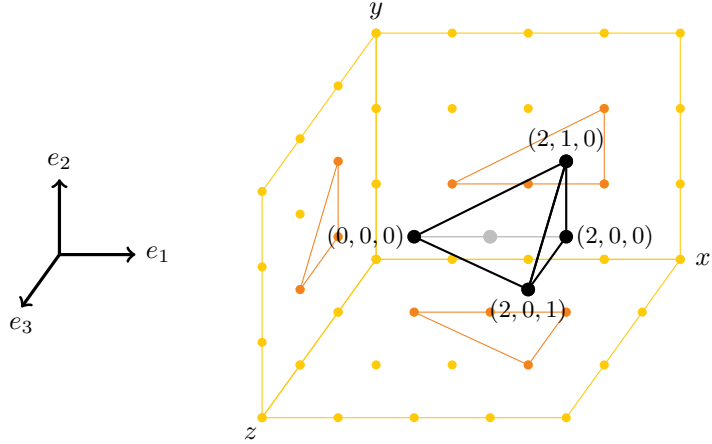}
	\caption{An example explaining how polytopes in $\R^3$ are visualized in our figures. The $x$-axis is directed to the right, the $y$-axis is directed upwards and the $z$-axis is directed towards the observer. Each figure depicts appropriately chosen planes orthogonal to the coordinate axes and the orthogonal projections of the polytope onto these planes. In this figure, the planes are given by equations $x=-1$, $y=-1$ and $z=-1$. \label{fig:expl:illustration}}
\end{figure}

We present each polytope as the Minkowski sum of indecomposable lattice polytopes. These
representations of are obtained using Magma~\cite{magma}.
\end{adjustwidth}

\bigskip

\mycaption{Volume 1}

\begin{longtable}{C{0.5em} C{15em} C{2em} C{0.5em} C{15em}} 
1. & \begin{tabular}{C{15em}} 
\begin{tikzpicture}[scale=0.5,baseline=-5pt] 
\draw[yellow!80!red] (-0.500000000000000, -0.300000000000000) -- (-0.500000000000000, 2.70000000000000);\draw[yellow!80!red] (-0.500000000000000, -0.300000000000000) -- (2.50000000000000, -0.300000000000000);\draw[yellow!80!red] (-0.500000000000000, 2.70000000000000) -- (2.50000000000000, 2.70000000000000);\draw[yellow!80!red] (2.50000000000000, -0.300000000000000) -- (2.50000000000000, 2.70000000000000);\fill[yellow!80!red] (-0.500000000000000, -0.300000000000000) circle (0.060000);\fill[yellow!80!red] (-0.500000000000000, 0.700000000000000) circle (0.060000);\fill[yellow!80!red] (-0.500000000000000, 1.70000000000000) circle (0.060000);\fill[yellow!80!red] (-0.500000000000000, 2.70000000000000) circle (0.060000);\fill[yellow!80!red] (0.500000000000000, -0.300000000000000) circle (0.060000);\fill[yellow!80!red] (0.500000000000000, 0.700000000000000) circle (0.060000);\fill[yellow!80!red] (0.500000000000000, 1.70000000000000) circle (0.060000);\fill[yellow!80!red] (0.500000000000000, 2.70000000000000) circle (0.060000);\fill[yellow!80!red] (1.50000000000000, -0.300000000000000) circle (0.060000);\fill[yellow!80!red] (1.50000000000000, 0.700000000000000) circle (0.060000);\fill[yellow!80!red] (1.50000000000000, 1.70000000000000) circle (0.060000);\fill[yellow!80!red] (1.50000000000000, 2.70000000000000) circle (0.060000);\fill[yellow!80!red] (2.50000000000000, -0.300000000000000) circle (0.060000);\fill[yellow!80!red] (2.50000000000000, 0.700000000000000) circle (0.060000);\fill[yellow!80!red] (2.50000000000000, 1.70000000000000) circle (0.060000);\fill[yellow!80!red] (2.50000000000000, 2.70000000000000) circle (0.060000);\draw[yellow!80!red] (-0.500000000000000, -0.300000000000000) -- (-2.00000000000000, -2.40000000000000);\draw[yellow!80!red] (-0.500000000000000, -0.300000000000000) -- (2.50000000000000, -0.300000000000000);\draw[yellow!80!red] (-2.00000000000000, -2.40000000000000) -- (1.00000000000000, -2.40000000000000);\draw[yellow!80!red] (2.50000000000000, -0.300000000000000) -- (1.00000000000000, -2.40000000000000);\fill[yellow!80!red] (-0.500000000000000, -0.300000000000000) circle (0.060000);\fill[yellow!80!red] (-1.00000000000000, -1.00000000000000) circle (0.060000);\fill[yellow!80!red] (-1.50000000000000, -1.70000000000000) circle (0.060000);\fill[yellow!80!red] (-2.00000000000000, -2.40000000000000) circle (0.060000);\fill[yellow!80!red] (0.500000000000000, -0.300000000000000) circle (0.060000);\fill[yellow!80!red] (0.000000000000000, -1.00000000000000) circle (0.060000);\fill[yellow!80!red] (-0.500000000000000, -1.70000000000000) circle (0.060000);\fill[yellow!80!red] (-1.00000000000000, -2.40000000000000) circle (0.060000);\fill[yellow!80!red] (1.50000000000000, -0.300000000000000) circle (0.060000);\fill[yellow!80!red] (1.00000000000000, -1.00000000000000) circle (0.060000);\fill[yellow!80!red] (0.500000000000000, -1.70000000000000) circle (0.060000);\fill[yellow!80!red] (0.000000000000000, -2.40000000000000) circle (0.060000);\fill[yellow!80!red] (2.50000000000000, -0.300000000000000) circle (0.060000);\fill[yellow!80!red] (2.00000000000000, -1.00000000000000) circle (0.060000);\fill[yellow!80!red] (1.50000000000000, -1.70000000000000) circle (0.060000);\fill[yellow!80!red] (1.00000000000000, -2.40000000000000) circle (0.060000);\draw[yellow!80!red] (-0.500000000000000, -0.300000000000000) -- (-2.00000000000000, -2.40000000000000);\draw[yellow!80!red] (-0.500000000000000, -0.300000000000000) -- (-0.500000000000000, 2.70000000000000);\draw[yellow!80!red] (-2.00000000000000, -2.40000000000000) -- (-2.00000000000000, 0.600000000000000);\draw[yellow!80!red] (-0.500000000000000, 2.70000000000000) -- (-2.00000000000000, 0.600000000000000);\fill[yellow!80!red] (-0.500000000000000, -0.300000000000000) circle (0.060000);\fill[yellow!80!red] (-1.00000000000000, -1.00000000000000) circle (0.060000);\fill[yellow!80!red] (-1.50000000000000, -1.70000000000000) circle (0.060000);\fill[yellow!80!red] (-2.00000000000000, -2.40000000000000) circle (0.060000);\fill[yellow!80!red] (-0.500000000000000, 0.700000000000000) circle (0.060000);\fill[yellow!80!red] (-1.00000000000000, 0.000000000000000) circle (0.060000);\fill[yellow!80!red] (-1.50000000000000, -0.700000000000000) circle (0.060000);\fill[yellow!80!red] (-2.00000000000000, -1.40000000000000) circle (0.060000);\fill[yellow!80!red] (-0.500000000000000, 1.70000000000000) circle (0.060000);\fill[yellow!80!red] (-1.00000000000000, 1.00000000000000) circle (0.060000);\fill[yellow!80!red] (-1.50000000000000, 0.300000000000000) circle (0.060000);\fill[yellow!80!red] (-2.00000000000000, -0.400000000000000) circle (0.060000);\fill[yellow!80!red] (-0.500000000000000, 2.70000000000000) circle (0.060000);\fill[yellow!80!red] (-1.00000000000000, 2.00000000000000) circle (0.060000);\fill[yellow!80!red] (-1.50000000000000, 1.30000000000000) circle (0.060000);\fill[yellow!80!red] (-2.00000000000000, 0.600000000000000) circle (0.060000);\node at (2.80000000000000, -0.300000000000000) {$x$};\node at (-0.500000000000000, 3.00000000000000) {$y$};\node at (-2.15000000000000, -2.61000000000000) {$z$};\draw[yellow!40!red] (0.000000000000000, -1.00000000000000) -- (-0.500000000000000, -1.70000000000000);\draw[yellow!40!red] (0.000000000000000, -1.00000000000000) -- (1.00000000000000, -1.00000000000000);\draw[yellow!40!red] (-0.500000000000000, -1.70000000000000) -- (1.00000000000000, -1.00000000000000);\fill[yellow!40!red] (0.000000000000000, -1.00000000000000) circle (0.060000);\fill[yellow!40!red] (-0.500000000000000, -1.70000000000000) circle (0.060000);\fill[yellow!40!red] (1.00000000000000, -1.00000000000000) circle (0.060000);\draw[yellow!40!red] (-1.00000000000000, 0.000000000000000) -- (-1.50000000000000, -0.700000000000000);\draw[yellow!40!red] (-1.00000000000000, 0.000000000000000) -- (-1.00000000000000, 1.00000000000000);\draw[yellow!40!red] (-1.50000000000000, -0.700000000000000) -- (-1.00000000000000, 1.00000000000000);\fill[yellow!40!red] (-1.00000000000000, 0.000000000000000) circle (0.060000);\fill[yellow!40!red] (-1.50000000000000, -0.700000000000000) circle (0.060000);\fill[yellow!40!red] (-1.00000000000000, 1.00000000000000) circle (0.060000);\draw[yellow!40!red] (0.500000000000000, 0.700000000000000) -- (0.500000000000000, 1.70000000000000);\draw[yellow!40!red] (0.500000000000000, 0.700000000000000) -- (1.50000000000000, 0.700000000000000);\draw[yellow!40!red] (0.500000000000000, 1.70000000000000) -- (1.50000000000000, 0.700000000000000);\fill[yellow!40!red] (0.500000000000000, 0.700000000000000) circle (0.060000);\fill[yellow!40!red] (0.500000000000000, 1.70000000000000) circle (0.060000);\fill[yellow!40!red] (1.50000000000000, 0.700000000000000) circle (0.060000);\draw[lightgray] (0.000000000000000, 0.000000000000000) -- (-0.500000000000000, -0.700000000000000);\draw[lightgray] (0.000000000000000, 0.000000000000000) -- (0.000000000000000, 1.00000000000000);\draw[lightgray] (-0.500000000000000, -0.700000000000000) -- (0.000000000000000, 1.00000000000000);\fill[lightgray] (0.000000000000000, 0.000000000000000) circle (0.090000);\fill[lightgray] (-0.500000000000000, -0.700000000000000) circle (0.090000);\fill[lightgray] (0.000000000000000, 1.00000000000000) circle (0.090000);\draw[lightgray] (0.000000000000000, 0.000000000000000) -- (-0.500000000000000, -0.700000000000000);\draw[lightgray] (0.000000000000000, 0.000000000000000) -- (1.00000000000000, 0.000000000000000);\draw[lightgray] (-0.500000000000000, -0.700000000000000) -- (1.00000000000000, 0.000000000000000);\fill[lightgray] (0.000000000000000, 0.000000000000000) circle (0.090000);\fill[lightgray] (-0.500000000000000, -0.700000000000000) circle (0.090000);\fill[lightgray] (1.00000000000000, 0.000000000000000) circle (0.090000);\draw[lightgray] (0.000000000000000, 0.000000000000000) -- (0.000000000000000, 1.00000000000000);\draw[lightgray] (0.000000000000000, 0.000000000000000) -- (1.00000000000000, 0.000000000000000);\draw[lightgray] (0.000000000000000, 1.00000000000000) -- (1.00000000000000, 0.000000000000000);\fill[lightgray] (0.000000000000000, 0.000000000000000) circle (0.090000);\fill[lightgray] (0.000000000000000, 1.00000000000000) circle (0.090000);\fill[lightgray] (1.00000000000000, 0.000000000000000) circle (0.090000);\draw[black,thick] (-0.500000000000000, -0.700000000000000) -- (0.000000000000000, 1.00000000000000);\draw[black,thick] (-0.500000000000000, -0.700000000000000) -- (1.00000000000000, 0.000000000000000);\draw[black,thick] (0.000000000000000, 1.00000000000000) -- (1.00000000000000, 0.000000000000000);\fill[black] (-0.500000000000000, -0.700000000000000) circle (0.090000);\fill[black] (0.000000000000000, 1.00000000000000) circle (0.090000);\fill[black] (1.00000000000000, 0.000000000000000) circle (0.090000);
\end{tikzpicture}\\ 
$\conv(0, e_1, e_2, e_3)$\end{tabular}\end{longtable} 

\mycaption{Volume 2}

\input{dim_3_mv_2_single_polytopes.tex}

\mycaption{Volume 3}

\input{dim_3_mv_3_rmax_single_polytopes.tex}

\mycaption{Volume 4}

\input{dim_3_mv_4_rmax_single_polytopes.tex}

\mycaption{Mixed Volume 2, full-dimensional, $\R$-maximal}

\input{dim_3_mv_2.tex}

\newpage
\mycaption{Mixed Volume 2, lower-dimensional, $\R$-maximal}

\input{dim_3_mv_2_lower_d_rmax.tex}

\mycaption{Mixed Volume 2, lower-dimensional, $\Z$-maximal but not $\R$-maximal}

\input{dim_3_mv_2_lower_d_zmax.tex}

\mycaption{Mixed Volume 3, full-dimensional, $\R$-maximal}

\input{dim_3_mv_3_rmax.tex}

\mycaption{Mixed Volume 3, full-dimensional, $\Z$-maximal but not $\R$-maximal}

\input{dim_3_mv_3_zmax.tex}

\mycaption{Mixed Volume 3, lower-dimensional, $\R$-maximal}
\input{dim_3_mv_3_lower_d_rmax.tex}

\mycaption{Mixed Volume 3, lower-dimensional, $\Z$-maximal but not $\R$-maximal}
\input{dim_3_mv_3_lower_d_zmax.tex}

\mycaption{Mixed Volume 4, full-dimensional, $\R$-maximal}

\input{dim_3_mv_4_rmax.tex}

\mycaption{Mixed Volume 4, full-dimensional, $\Z$-maximal but not $\R$-maximal}

 \input{dim_3_mv_4_zmax.tex}

\newpage
\mycaption{Mixed Volume 4, lower-dimensional, $\R$-maximal}

\input{dim_3_mv_4_lower_d_rmax.tex}

\mycaption{Mixed Volume 4, lower-dimensional, $\Z$-maximal but not $\R$-maximal}

\input{dim_3_mv_4_lower_d_zmax.tex}

\section{Enumeration Data for dimension 2}\label{appendix:b}

\begin{adjustwidth}{2cm}{2cm}
For the reader's convenience we present a complete list 
of maximal pairs of lattice polygons of mixed volume up to 4.
Note that these have already been classified in \cite{EsterovGusev2016}.
\end{adjustwidth}

\bigskip

\mycaption{Mixed Volume 1}

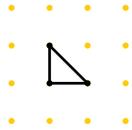
\begin{longtable}{C{0.5em} C{15em} C{2em} C{15em}} 
1. & \begin{tikzpicture}[scale=0.5,baseline=-5pt] 
\fill[yellow!80!red] (-1,-1) circle (0.08);\fill[yellow!80!red] (-1,0) circle (0.08);\fill[yellow!80!red] (-1,1) circle (0.08);\fill[yellow!80!red] (-1,2) circle (0.08);\fill[yellow!80!red] (0,-1) circle (0.08);\fill[yellow!80!red] (0,0) circle (0.08);\fill[yellow!80!red] (0,1) circle (0.08);\fill[yellow!80!red] (0,2) circle (0.08);\fill[yellow!80!red] (1,-1) circle (0.08);\fill[yellow!80!red] (1,0) circle (0.08);\fill[yellow!80!red] (1,1) circle (0.08);\fill[yellow!80!red] (1,2) circle (0.08);\fill[yellow!80!red] (2,-1) circle (0.08);\fill[yellow!80!red] (2,0) circle (0.08);\fill[yellow!80!red] (2,1) circle (0.08);\fill[yellow!80!red] (2,2) circle (0.08); 
\draw[line width=1pt] (0, 0) -- (0, 1);\draw[line width=1pt] (0, 0) -- (1, 0);\draw[line width=1pt] (0, 1) -- (1, 0);\node at (1/3, 1/3) {}; 
\fill (0,0) circle (0.08);\fill (0,1) circle (0.08);\fill (1,0) circle (0.08);\end{tikzpicture} & & \begin{tikzpicture}[scale=0.5,baseline=-5pt] 
\fill[yellow!80!red] (-1,-1) circle (0.08);\fill[yellow!80!red] (-1,0) circle (0.08);\fill[yellow!80!red] (-1,1) circle (0.08);\fill[yellow!80!red] (-1,2) circle (0.08);\fill[yellow!80!red] (0,-1) circle (0.08);\fill[yellow!80!red] (0,0) circle (0.08);\fill[yellow!80!red] (0,1) circle (0.08);\fill[yellow!80!red] (0,2) circle (0.08);\fill[yellow!80!red] (1,-1) circle (0.08);\fill[yellow!80!red] (1,0) circle (0.08);\fill[yellow!80!red] (1,1) circle (0.08);\fill[yellow!80!red] (1,2) circle (0.08);\fill[yellow!80!red] (2,-1) circle (0.08);\fill[yellow!80!red] (2,0) circle (0.08);\fill[yellow!80!red] (2,1) circle (0.08);\fill[yellow!80!red] (2,2) circle (0.08); 
\draw[line width=1pt] (0, 0) -- (0, 1);\draw[line width=1pt] (0, 0) -- (1, 0);\draw[line width=1pt] (0, 1) -- (1, 0);\node at (1/3, 1/3) {}; 
\fill (0,0) circle (0.08);\fill (0,1) circle (0.08);\fill (1,0) circle (0.08);\end{tikzpicture}\\ \\ 
 & $\conv(0, e_1, e_2)$ & & $\conv(0, e_1, e_2)$\\ \\ \\ \\ 
\end{longtable} 

\mycaption{Mixed Volume 2}

\begin{longtable}{C{0.5em} C{15em} C{2em} C{15em}} 
1. & \begin{tikzpicture}[scale=0.5,baseline=-5pt] 
\fill[yellow!80!red] (-1,-1) circle (0.08);\fill[yellow!80!red] (-1,0) circle (0.08);\fill[yellow!80!red] (-1,1) circle (0.08);\fill[yellow!80!red] (-1,2) circle (0.08);\fill[yellow!80!red] (-1,3) circle (0.08);\fill[yellow!80!red] (0,-1) circle (0.08);\fill[yellow!80!red] (0,0) circle (0.08);\fill[yellow!80!red] (0,1) circle (0.08);\fill[yellow!80!red] (0,2) circle (0.08);\fill[yellow!80!red] (0,3) circle (0.08);\fill[yellow!80!red] (1,-1) circle (0.08);\fill[yellow!80!red] (1,0) circle (0.08);\fill[yellow!80!red] (1,1) circle (0.08);\fill[yellow!80!red] (1,2) circle (0.08);\fill[yellow!80!red] (1,3) circle (0.08);\fill[yellow!80!red] (2,-1) circle (0.08);\fill[yellow!80!red] (2,0) circle (0.08);\fill[yellow!80!red] (2,1) circle (0.08);\fill[yellow!80!red] (2,2) circle (0.08);\fill[yellow!80!red] (2,3) circle (0.08); 
\draw[line width=1pt] (0, 0) -- (0, 2);\draw[line width=1pt] (0, 0) -- (1, 0);\draw[line width=1pt] (0, 2) -- (1, 0);\node at (1/3, 2/3) {}; 
\fill (0,0) circle (0.08);\fill (0,1) circle (0.08);\fill (0,2) circle (0.08);\fill (1,0) circle (0.08);\end{tikzpicture} & & \begin{tikzpicture}[scale=0.5,baseline=-5pt] 
\fill[yellow!80!red] (-1,-1) circle (0.08);\fill[yellow!80!red] (-1,0) circle (0.08);\fill[yellow!80!red] (-1,1) circle (0.08);\fill[yellow!80!red] (-1,2) circle (0.08);\fill[yellow!80!red] (-1,3) circle (0.08);\fill[yellow!80!red] (0,-1) circle (0.08);\fill[yellow!80!red] (0,0) circle (0.08);\fill[yellow!80!red] (0,1) circle (0.08);\fill[yellow!80!red] (0,2) circle (0.08);\fill[yellow!80!red] (0,3) circle (0.08);\fill[yellow!80!red] (1,-1) circle (0.08);\fill[yellow!80!red] (1,0) circle (0.08);\fill[yellow!80!red] (1,1) circle (0.08);\fill[yellow!80!red] (1,2) circle (0.08);\fill[yellow!80!red] (1,3) circle (0.08);\fill[yellow!80!red] (2,-1) circle (0.08);\fill[yellow!80!red] (2,0) circle (0.08);\fill[yellow!80!red] (2,1) circle (0.08);\fill[yellow!80!red] (2,2) circle (0.08);\fill[yellow!80!red] (2,3) circle (0.08); 
\draw[line width=1pt] (0, 0) -- (0, 2);\draw[line width=1pt] (0, 0) -- (1, 0);\draw[line width=1pt] (0, 2) -- (1, 0);\node at (1/3, 2/3) {}; 
\fill (0,0) circle (0.08);\fill (0,1) circle (0.08);\fill (0,2) circle (0.08);\fill (1,0) circle (0.08);\end{tikzpicture}\\ \\ 
 & $\conv(0, e_1,2 e_2)$ & & $\conv(0, e_1,2 e_2)$\\ \\ \\ \\ 
2. & \begin{tikzpicture}[scale=0.5,baseline=-5pt] 
\fill[yellow!80!red] (-1,-2) circle (0.08);\fill[yellow!80!red] (-1,-1) circle (0.08);\fill[yellow!80!red] (-1,0) circle (0.08);\fill[yellow!80!red] (-1,1) circle (0.08);\fill[yellow!80!red] (-1,2) circle (0.08);\fill[yellow!80!red] (0,-2) circle (0.08);\fill[yellow!80!red] (0,-1) circle (0.08);\fill[yellow!80!red] (0,0) circle (0.08);\fill[yellow!80!red] (0,1) circle (0.08);\fill[yellow!80!red] (0,2) circle (0.08);\fill[yellow!80!red] (1,-2) circle (0.08);\fill[yellow!80!red] (1,-1) circle (0.08);\fill[yellow!80!red] (1,0) circle (0.08);\fill[yellow!80!red] (1,1) circle (0.08);\fill[yellow!80!red] (1,2) circle (0.08);\fill[yellow!80!red] (2,-2) circle (0.08);\fill[yellow!80!red] (2,-1) circle (0.08);\fill[yellow!80!red] (2,0) circle (0.08);\fill[yellow!80!red] (2,1) circle (0.08);\fill[yellow!80!red] (2,2) circle (0.08); 
\draw[line width=1pt] (0, 0) -- (0, 1);\draw[line width=1pt] (0, 0) -- (1, -1);\draw[line width=1pt] (0, 1) -- (1, 0);\draw[line width=1pt] (1, -1) -- (1, 0);\node at (1/2, 0) {}; 
\fill (0,0) circle (0.08);\fill (0,1) circle (0.08);\fill (1,-1) circle (0.08);\fill (1,0) circle (0.08);\end{tikzpicture} & & \begin{tikzpicture}[scale=0.5,baseline=-5pt] 
\fill[yellow!80!red] (-1,-2) circle (0.08);\fill[yellow!80!red] (-1,-1) circle (0.08);\fill[yellow!80!red] (-1,0) circle (0.08);\fill[yellow!80!red] (-1,1) circle (0.08);\fill[yellow!80!red] (-1,2) circle (0.08);\fill[yellow!80!red] (0,-2) circle (0.08);\fill[yellow!80!red] (0,-1) circle (0.08);\fill[yellow!80!red] (0,0) circle (0.08);\fill[yellow!80!red] (0,1) circle (0.08);\fill[yellow!80!red] (0,2) circle (0.08);\fill[yellow!80!red] (1,-2) circle (0.08);\fill[yellow!80!red] (1,-1) circle (0.08);\fill[yellow!80!red] (1,0) circle (0.08);\fill[yellow!80!red] (1,1) circle (0.08);\fill[yellow!80!red] (1,2) circle (0.08);\fill[yellow!80!red] (2,-2) circle (0.08);\fill[yellow!80!red] (2,-1) circle (0.08);\fill[yellow!80!red] (2,0) circle (0.08);\fill[yellow!80!red] (2,1) circle (0.08);\fill[yellow!80!red] (2,2) circle (0.08); 
\draw[line width=1pt] (0, 0) -- (0, 1);\draw[line width=1pt] (0, 0) -- (1, -1);\draw[line width=1pt] (0, 1) -- (1, 0);\draw[line width=1pt] (1, -1) -- (1, 0);\node at (1/2, 0) {}; 
\fill (0,0) circle (0.08);\fill (0,1) circle (0.08);\fill (1,-1) circle (0.08);\fill (1,0) circle (0.08);\end{tikzpicture}\\ \\ 
 & $\conv(0, e_1- e_2)+\conv(0, e_2)$ & & $\conv(0, e_1- e_2)+\conv(0, e_2)$\\ \\ \\ \\ 
3. & \begin{tikzpicture}[scale=0.5,baseline=-5pt] 
\fill[yellow!80!red] (-1,-1) circle (0.08);\fill[yellow!80!red] (-1,0) circle (0.08);\fill[yellow!80!red] (-1,1) circle (0.08);\fill[yellow!80!red] (-1,2) circle (0.08);\fill[yellow!80!red] (0,-1) circle (0.08);\fill[yellow!80!red] (0,0) circle (0.08);\fill[yellow!80!red] (0,1) circle (0.08);\fill[yellow!80!red] (0,2) circle (0.08);\fill[yellow!80!red] (1,-1) circle (0.08);\fill[yellow!80!red] (1,0) circle (0.08);\fill[yellow!80!red] (1,1) circle (0.08);\fill[yellow!80!red] (1,2) circle (0.08);\fill[yellow!80!red] (2,-1) circle (0.08);\fill[yellow!80!red] (2,0) circle (0.08);\fill[yellow!80!red] (2,1) circle (0.08);\fill[yellow!80!red] (2,2) circle (0.08); 
\draw[line width=1pt] (0, 0) -- (0, 1);\draw[line width=1pt] (0, 0) -- (1, 0);\draw[line width=1pt] (0, 1) -- (1, 0);\node at (1/3, 1/3) {}; 
\fill (0,0) circle (0.08);\fill (0,1) circle (0.08);\fill (1,0) circle (0.08);\end{tikzpicture} & & \begin{tikzpicture}[scale=0.5,baseline=-5pt] 
\fill[yellow!80!red] (-1,-1) circle (0.08);\fill[yellow!80!red] (-1,0) circle (0.08);\fill[yellow!80!red] (-1,1) circle (0.08);\fill[yellow!80!red] (-1,2) circle (0.08);\fill[yellow!80!red] (-1,3) circle (0.08);\fill[yellow!80!red] (0,-1) circle (0.08);\fill[yellow!80!red] (0,0) circle (0.08);\fill[yellow!80!red] (0,1) circle (0.08);\fill[yellow!80!red] (0,2) circle (0.08);\fill[yellow!80!red] (0,3) circle (0.08);\fill[yellow!80!red] (1,-1) circle (0.08);\fill[yellow!80!red] (1,0) circle (0.08);\fill[yellow!80!red] (1,1) circle (0.08);\fill[yellow!80!red] (1,2) circle (0.08);\fill[yellow!80!red] (1,3) circle (0.08);\fill[yellow!80!red] (2,-1) circle (0.08);\fill[yellow!80!red] (2,0) circle (0.08);\fill[yellow!80!red] (2,1) circle (0.08);\fill[yellow!80!red] (2,2) circle (0.08);\fill[yellow!80!red] (2,3) circle (0.08);\fill[yellow!80!red] (3,-1) circle (0.08);\fill[yellow!80!red] (3,0) circle (0.08);\fill[yellow!80!red] (3,1) circle (0.08);\fill[yellow!80!red] (3,2) circle (0.08);\fill[yellow!80!red] (3,3) circle (0.08); 
\draw[line width=1pt] (0, 0) -- (0, 2);\draw[line width=1pt] (0, 0) -- (2, 0);\draw[line width=1pt] (0, 2) -- (2, 0);\node at (2/3, 2/3) {}; 
\fill (0,0) circle (0.08);\fill (0,1) circle (0.08);\fill (0,2) circle (0.08);\fill (1,0) circle (0.08);\fill (1,1) circle (0.08);\fill (2,0) circle (0.08);\end{tikzpicture}\\ \\ 
 & $\conv(0, e_1, e_2)$ & & $2\conv(0, e_1, e_2)$\\ \\ \\ \\ 
\end{longtable} 

\mycaption{Mixed Volume 3}

\input{dim_2_mv_3.tex}

\mycaption{Mixed Volume 4}

\input{dim_2_mv_4.tex}

\newpage
\restoregeometry

\end{document}